\documentclass[12pt]{amsart}
\newtheorem{theo}{{\bfseries Theorem}}[section]
\newtheorem{prop}[theo]{{\bfseries Proposition}}
\newtheorem{lem}[theo]{{\bfseries Lemma}}
\newtheorem{cor}[theo]{{\bfseries Corollary}}
\newtheorem{df}[theo]{{\bfseries Definition}}

\input xy

\def \N {\mathbb N}

\def \a {\alpha }
\def \b {\beta}

\def \ep {\epsilon}

\def \t {\theta}

\usepackage{amssymb,latexsym}
\usepackage[mathcal]{eucal}
\usepackage{amscd}
\usepackage{amsmath}
\usepackage{graphicx}
\usepackage{amsfonts}
\usepackage{amscd}

\numberwithin{equation}{section}

\makeindex

\begin{document}

\title{\bfseries  Maximal r-Diameter Sets and \\ Solids of Constant Width}

\author{Ethan Akin}

\address{Mathematics Department \\
    The City College \\
    137 Street and Convent Avenue \\
    New York City, NY 10031, USA}

\date{January, 2010}
\vspace{.5cm} \maketitle

{\bfseries Abstract:}  We recall the definition of an $r$-maximal set in a metric space as a maximal subset of diameter $r$.
In the special case when the metric space is Euclidean such a set is exactly a solid of constant diameter $r$.
In the process of reviewing the theory of these objects we provide a simple construction which generates a large class of such solids.
\vspace{.5cm}


\tableofcontents

\section*{Introduction}

I was reading a paper - reviewing it actually - in which a subset was being covered by sets of diameter at most $r$ for
some fixed $r > 0$.
It occurred to me that if a subset is rather elongated then the diameter constraint is not completely binding in the
sense that we can enlarge the set laterally without increasing the diameter.
This suggests the question, which can be asked
in any metric space $X$, of characterizing the subsets of diameter $r$ which are maximal with respect to this condition.

Eventually I discovered that for Euclidean space this notion is quite old.  Under the name \emph{completeness} or
\emph{diametrical completeness} the concept was introduced by Meissner in 1909 and its relation with the
related concept of \emph{constant width} has been the object of considerable study.  Most of what I rediscovered
had been analyzed by Minkowski and was described in Bonnesen and Fenchel's (1934) survey, see also Eggleston (1958) and
Lay (1987).  We will follow Eggleston (1965) in using  \emph{diametrical maximality} in place of the overused
term \emph{completeness}.  We also overlap recent work in Lachand-Robert and Oudet (2007)
and Bayen, Lachand-Robert and Oudet (2007).

From Zorn's Lemma it follows that any set of diameter at most $r$, which we will call an \emph{$r$-bounded} set, is a
subset of some maximal subset of diameter $r$.  We will call such a set an \emph{$r$-maximal} set.

In an $r$-bounded set $C$ we will call a pair $x_1, x_2 \in C$ \emph{antipodal} if $d(x_1,x_2) = r$.

Conditions on $r$-bounded sets depend on the underlying metric.  Call the metric $d$  \emph{open}(or \emph{proper}) when
for all $x \in X$ the function from $X$ to $[0,\infty)$ defined by $y \mapsto d(x,y)$ is an open map (resp. a proper
map).  Thus, $d$ is proper iff any closed, finite diameter subset is compact. For any compact
subset $C$ of $X$ define $d_+(C,x) = max \{ d(x,y) : y \in C \}$. Equivalently, $d_+(C,x)$ is the radius
of the smallest closed ball centered at $x$ which contains $C$.   Call $d$ \emph{connected} when
for every compact $C$ and any $r > 0$ the open set $\{ x \in X : d_+(C,x) < r \}$ is connected and, when it is
nonempty,  its closure is $\{ x \in X : d_+(C,x) \leq r \}$.

\begin{theo} Let $X$ be a metric space with a proper, open  metric.

(a) If  $C$ is an $r$-maximal subset then $C$ is compact and for every
$x_1$ in the topological boundary of $C$ there is an antipodal point $x_2$ in the boundary of $C$.  Furthermore,
$C = \{ x \in X : d_+(C,x) \leq r \}$ and the interior of $C$ is $\{ x \in X : d_+(C,x) < r \}$. If the
metric is connected then $C$ is connected.

(b) Assume that the metric is  connected and that $C$ is a compact $r$-bounded subset with
nonempty interior. If every point of the boundary of $C$ has an antipode in $C$ then $C$ is $r$-maximal.
\end{theo}
\vspace{.5cm}

If $E$ is a Minkowski space, i.e. a finite dimensional, normed linear space then the
associated metric is proper, open and connected.
For a functional $\omega$ in the unit sphere $S^*$ of the dual space of $E$ we define the $\omega$ diameter
of a subset $C$ to be $sup \{ \omega(x_1 - x_2) : x_1, x_2 \in C \}$. We say that $C$ has \emph{constant
diameter $r$} if the $\omega$ diameter of $C$ equals $r$ for every $\omega \in S^*$.

\begin{theo} Let $E$ be a Minkowski space.

(a) If $C$ is an $r$-maximal subset  then $C$ is a
compact convex set with nonempty interior.

(b) If $C$ is a compact, convex set of constant diameter $r$ then $C$ is an $r$-maximal subset.
\end{theo}

By a Euclidean space we mean a finite dimensional, linear space with norm obtained from an inner product.
When the dimension is $n$ such a space is isometric to ${\mathbb R}^n$ with the usual metric.

\begin{theo} Let $E$ be a Euclidean space.

(a) A compact convex subset $C$ is $r$-maximal iff it has constant diameter $r$.

(b) If a point $x$ of an $r$-maximal set $C$ has two distinct antipodal points $y_1, y_2$ in $C$ then the
arc between them centered at $x$ is entirely contained in the boundary of $C$.  Conversely, if an
arc of radius $r$ is contained in the boundary of an $r$-maximal set $C$ then the center $x$ lies in the
boundary of $C$ and for each point $y$ of the arc other than the endpoints $x$ is the unique point of $C$
antipodal to $y$.
\end{theo}

In the Euclidean  case there is a simple construction which yields a large class of examples and which appears to be new.

Let $g $ be a $C^2$ real-valued function on the unit sphere $S$ in ${\mathbb R}^n$ ($n \geq 2$). Assume that
$g$ is an odd function, i.e. $g(-u) = - g(u)$ for all $u \in S$.  Extend $g$ to
a $C^2$ odd, homogeneous function of degree $1$ defined on ${\mathbb R}^n \setminus 0$ by $g(x) = |x|g(x/|x|)$.
Let $H = (h_1,...,h_n) : {\mathbb R}^n \setminus 0 \to {\mathbb R}^n $ be the gradient of $g$ so
that $h_i = \partial g/\partial x_i$
for $i = 1,..,n$. Each $h_i$ is a $C^1$ even, homogeneous function of degree $0$. We say that $g$ satisfies the
\emph{$r$-Maximality Condition} when the Jacobian matrix $( \partial^2 g/\partial x_i \partial x_j)$ has all
eigenvalues contained in $[-(r/2),(r/2)]$ at every point $u$ of $S$.
Notice that for $g$ any $C^2$ odd function
on $S$ there is a maximum positive $\lambda^*$ such that $\lambda^* g$ satisfies the $r$-Rotundity Condition
and then $\lambda g$ satisfies  $r$-Maximality Condition for all $\lambda \in [0, \lambda^*]$.

\begin{theo}  If $g$ satisfies the $r$-Maximality  Condition then the subset of ${\mathbb R}^n$
\begin{displaymath}
\begin{split}
C  \ =  \  \{ (H(u) + t (r/2)u) : (u,t) \in S \times [-1,1]  \} \\ =
\  \{ (H(u) + t (r/2) u) : (u,t) \in S \times [0,1]  \} \hspace{.5cm}
\end{split}\end{displaymath}
is a solid of constant diameter $r$, and so is an $r$-maximal set.
\end{theo}
\vspace{.5cm}

In general, if $C$ is a solid of constant width $r$ then for each $u \in S$
there is a unique directed segment of length $r$ with direction $u$
connecting a pair of antipodal points in $C$.  We can define $H(u)$ is the midpoint of the segment.
This is called the \emph{median surface} by Bayen, Lachand-Robert and Oudet (2007), see also Guilfoyle and Klingenberg
(2009). The function
  $H$ is even, i.e. $H(-u) = H(u)$, because the segment for $-u$ connects the same pair but with the
orientation reversed.  This is explains why the two descriptions of $C$ above yield the same set.

 Letting $\lambda$ vary in $[0,\lambda^*]$ we obtain a one parameter family of
$r$-maximal sets connecting $C$ to the unit ball centered at $0$.

In the planar case, a related construction shows that
the radius of curvature is the only constraint on including a piece of a curve in
the boundary of some $r$-maximal set.

\begin{theo}Given a plane curve with  radius of curvature bounded by $r$  at every point then any sufficiently
  short piece can be embedded in the boundary of some $r$-maximal subset of the plane. \end{theo}
\vspace{1cm}

\section{General Metric Spaces}

We recall the language of relations, see e.g. Akin (1993). For sets $X$ and $Y$
 a \emph{relation} $R$ from $X$ to $Y$ is an arbitrary subset $R \subset X \times Y$.
 For $x \in X, \ R(x) = \{ y \in Y : (x,y) \in R \}$ and for $A \subset X, \
 R(A) = \bigcup_{x \in A} \ R(x).$ The inverse relation  $R^{-1} $
from $Y$ to $X$ is $\{ (y,x) : (x,y) \in R \}$. The relation $R$ is
called \emph{surjective} when $R(X) = Y$ and $R^{-1}(Y) = X$.  That is, $R$ projects onto
each factor. When $X = Y$ we call $R \subset X \times X$  a relation on $X$.
If $X$ and $Y$ are topological spaces then $R$ is a closed (or open)
relation when it is a closed (resp. open) subset of the product.

Now let $X$  be a metric space with metric $d$.
With $r > 0$ define  relations on $X$
\begin{equation}\label{1.1}
\begin{split}
V_r = \{ (x,y) \in X \times X : d(x,y) < r \}, \\
\bar V_r = \{ (x,y) \in X \times X : d(x,y) \leq r \}, \\
A_r = \{ (x,y) \in X \times X : d(x,y) = r \}.
\end{split}
\end{equation}
For $C \subset X$ the diameter $diam(C) = sup \{ d(x,y) : x,y \in C \}$.
Thus, $diam (C) \leq r$ iff $C \times C \subset \bar V_r $.  Since $\bar V_r$ is closed it follows that
$diam (C) = diam (\overline{C})$ where $\overline{C}$ is the closure of $C$. $C$ is called \emph{bounded} when it
has finite diameter.

Let $Bdry(C)$ denote the
topological boundary of $C$ so that $Bdry(C) = \overline{C} \setminus C^{\circ}$, where $C^{\circ}$ is the
interior of $C$.

For a nonempty subset $C \subset X$ and a point $x \in X$ we define
\begin{equation}\label{1.1a}
\begin{split}
d_-(C,x) \quad =_{def} \quad  inf \{ d(y,x): y \in C \}, \hspace{2cm}\\
d_+(C,x) \quad =_{def} \quad  sup \{ d(y,x): y \in C \}. \hspace{2cm}
\end{split}
\end{equation}
Of course, $d_+(C,x)$ is only finite when $C$ is bounded.

For $C \subset D$ it is clear that
\begin{equation}\label{1.1b}
d_-(D,x) \ \leq \ d_-(C,x) \ \leq \ d_+(C,x) \ \leq \ d_+(D,x).
\end{equation}

For any subset $D \subset X$,
\begin{equation}\label{1.1c}
sup \{ d_-(C,x) : x \in D \} \quad = \quad inf \{ r \geq 0 : D \subset \bar V_r(C) \}.
\end{equation}
Of course, if $C$ and $D$ are compact then all of these sups and infs are achieved.

For compact subsets
$C,D$ the \emph{Hausdorff distance}  is defined to be
\begin{equation}\label{1.1d}
d(C,D) \quad =_{def} \quad min \{ r \geq 0 : D \subset \bar V_r(C) \ \mbox{and} \ C \subset \bar V_r(D) \}.
\end{equation}
For properties of the space of compact subsets equipped with the Hausdorff metric see Akin (1993) Chapter 7 or
Kuratowski (1968) Section 42.

From the triangle inequality it easily follows that for nonempty compacta $C, D \subset X$  and points $x, y \in X$
\begin{equation}\label{1.1e}
\begin{split}
|d_-(C,x) - d_-(D,y)| \quad \leq \quad d(C,D) + d(x,y); \hspace{1cm}\\
|d_+(C,x) - d_+(D,y)| \quad \leq \quad d(C,D) + d(x,y). \hspace{1cm}
\end{split}
\end{equation}

For any bounded subset $C$ of $X$ we define the \emph{antipodal relation} $A_C$ to
be the closed, symmetric relation on $C$
\begin{equation}\label{1.2}
A_C \quad = \quad A_r \cap (C \times C) \qquad \mbox{where} \qquad r \ = \ diam(C).
\end{equation}

Call the metric $d$ \emph{proper} when  for each $x \in X$ the
function $d(x, \cdot)$, i.e. $y \mapsto d(x,y)$ is a proper map from $X$ to $[0,\infty)$. That is,
 the pre-image of a compact set is
compact.  Equivalently, $d$ is proper when
any closed, bounded subset is compact.  Call $d$ \emph{open} when for each $x \in X$ the
function $ d(x,\cdot)$ is an open map from $X$ to $[0,\infty)$.  Notice that when the metric is proper, $X$ is
locally compact and if the metric is open then $X$ is not compact.

We can localize the rather strong condition that the metric be open. Call $r > 0$ a \emph{regular value} for
$x \in X$ when $r$ is neither a local maximum or a local minimum value for $d(x,\cdot)$. That is, if $d(x,y) = r$ then
$y$ is neither a local maximum point, nor a local minimum point for $d(x,\cdot)$.  Call $r$ a regular value for
$C \subset X$ when it is a regular value for every $x \in X$.  When the metric is open, $d(x,\cdot)$ has no
local maxima or positive local minima and so every positive $r$ is regular for every $x \in X$.

\begin{lem} \label{lem1.1} Let $C$ be a subset of $X$ with $diam(C) = r$ and let $x \in C$.

(a) If $x, y \in C$ and
$r$ is a regular value for $y$ then
  $d(x,y) = r$ implies $x \in Bdry(C)$.

(b) If $x \in C^{\circ}$ and $r$ is a regular value for $C$ then
 $A_C(x) = \emptyset$.
\end{lem}

{\bfseries Proof:} (a): Since $r$ is a regular value for $y$ there exist $z \in X$
arbitrarily close to $x$ with $r < d(z,y)$. Since $diam(C) = r, z \not\in C$.
Hence, $x \not\in C^{\circ}$.

(b): Apply (a) with $y \in C$ arbitrary.

$\Box$ \vspace{.5cm}

\begin{df} \label{df1.2}For $r > 0$ and a subset $C $ of $X$.  We will call $C$ \emph{$r$-bounded} when the
diameter of $C$ is at most $r$.  $C$ is called \emph{a diametrically maximal set of size $r$}, or just
 \emph{$r$-maximal}, when $C$ is a maximal $r$-bounded subset of $X$. \end{df}
\vspace{.5cm}

For $C$ a nonempty subset of $X$ and $r > 0$ define the $r$-dual
\begin{equation}\label{1.3}
\begin{split}
C^*_r \quad =_{def} \quad \{ y \in X : C \subset \bar V_r(y) \} \hspace{2cm} \\
 = \quad \bigcap \{ \bar V_r(x) : x \in C \} \quad = \quad \{ d_+(C,\cdot )  \leq r \}, \hspace{1cm}
\end{split}
\end{equation}
where we use the notation $\{ d_+(C,\cdot ) \leq r \}$ for $\{ y \in X : d_+(C,y) \leq r \}$ and similarly
for $\{ d_+(C,\cdot) < r \}$.

\begin{theo}\label{theo1.3} Let $C \subset X$.
\begin{itemize}
\item[(a)] For every $r > 0$, $C^*_r$ is a closed subset of $X$.
\item[(b)] If $C \subset D$ then $D^*_r \subset C^*_r$.
\item[(c)] $C \subset C^*_r$ if and only if $C$ is $r$-bounded.
\item[(d)] If $C$ is $r$-bounded then
\begin{equation}\label{1.4}
C^*_r \quad = \quad \{ y \in X : C \cup \{ y \} \ \mbox{is} \ r\mbox{-bounded} \  \}.
\end{equation}
\item[(e)]   $C$ is contained in an $r$-maximal set if and only if $C$ is $r$-bounded.
\item[(f)] $C$ is an $r$-maximal set if and only if $C = C^*_r$.
\item[(g)] If  $C$ is $r$-bounded then the following are equivalent:
\begin{enumerate}
\item[(i)] $C$ is contained in a unique $r$-maximal set.
\item[(ii)] $C^*_r$ is an $r$-maximal set.
\item[(iii)] $C^*_r$ is $r$-bounded.
\end{enumerate}
 When these conditions hold, $C^*_r$ is the unique $r$-maximal set containing $C$.
\end{itemize}
\end{theo}

{\bfseries Proof:} (a), (b) and (c): Obvious.

(d): If $diam(C) \leq r$ then $diam(C \cup \{ y \}) \leq r$ iff $d(y,x) \leq r$ for all $x \in C$ and so
iff $y \in C^*_r$.

(e): An r-round set is $r$-bounded and so any subset is $r$-bounded. On the other hand,
 the condition $C \times C \subset \bar V_r$ is a property of finite type. So it follows from Zorn's
Lemma that any $r$-bounded subset $C$  is contained in a maximal  $r$-bounded set.

(f): If $C$ is $r$-maximal or if $C = C^*_r$ then $C$ is $r$-bounded. Clearly, $C$ is $r$-maximal iff
$C =  \{ y \in X : diam(C \cup \{ y \}) \leq r \}$ and so iff $C = C^*_r$ by (d).

(g): Observe first that by (b) and (f)
\begin{equation}\label{1.5}
 C \subset D \quad \mbox{ and } \quad D  \ \mbox{is} \ r\mbox{-round} \qquad \Longrightarrow \qquad
D = D^*_r \subset C^*_r.
\end{equation}

(i) $\Rightarrow$ (ii): Let $D$ be the unique $r$-maximal set containing $C$ and let $y \in C^*_r$. By
 equation (\ref{1.4}) and (e) $\{y \} \cup C$ is contained in some $r$-maximal set which must be $D$ by
 uniqueness.  Hence, $y \in D$. As $y$ was arbitrary $C^*_r \subset D$. From implication (\ref{1.5})
 it follows that $C^*_r = D$ and so is itself $r$-maximal.

 (ii) $\Rightarrow$ (iii): An $r$-maximal set is $r$-bounded.

 (iii) $\Rightarrow$ (i): If $D$ is an $r$-maximal set which contains $C$ then by (\ref{1.5}) again
 $D \subset C^*_r$.  Since $C^*_r$ is $r$-bounded by assumption, maximality of $D$ implies $D = C^*_r$.
 Hence, $C^*_r$ is the unique $r$-maximal set which contains $C$.

$\Box$ \vspace{.5cm}

\begin{df} \label{def1.3a}Let $C$ be an $r$-bounded subset of $X$.  We say that
$C$ satisfies the \emph{antipodal condition}
if for every $x \in Bdry(C)$ there exists $y \in C$ such that $d(x,y) = r$, i.e. $y$ is antipodal to $x$.
Equivalently, $A_C(x) \not= \emptyset$ for
all $x \in Bdry(C)$. \end{df}
\vspace{.5cm}

Thus, if $r$ is a regular value for $C$, e.g. if the metric is open, then $C$ satisfies the antipodal condition
iff $A_C$ is a surjective relation on $Bdry(C)$.

\begin{prop}\label{prop1.4} If $C$ be an $r$-maximal subset of $X$ then $C$ is a closed subset of $X$.
 If $d$ is proper then $C$ is compact.
\end{prop}

{\bfseries Proof:}  $C = \overline{C}$ because taking closure does not increase diameter.  Alternatively,
apply Theorem \ref{theo1.3}(a) and (f). If $d$ is proper then $C$ is compact because it is closed and bounded.

$\Box$ \vspace{.5cm}

\begin{theo}\label{theo1.4a} Assume that $d$ is proper. If $\{C_n \}$ is a sequence of $r$-maximal sets which
converges to a compact set $C$ with respect to the Hausdorff metric, then $C$ satisfies the antipodal condition.
In particular, an $r$-maximal set satisfies the antipodal condition.\end{theo}

{\bfseries Proof:} Because the metric is proper and the sequence $\{ C_n \}$ is convergent, the
set $X_0 \ = \ \overline{\bigcup_n C_n}$ is closed and bounded and so is compact.  Furthermore,
the limit $C$ is the $limsup$ of the sequence.  That is,
\begin{equation}\label{1.5a}
C \quad  = \quad \bigcap_n \overline{ \bigcup_{k \geq n} \ C_k}.
\end{equation}
See Akin (1993) Lemma 7.5.

Fix $x \in Bdry(C)$ and let $\{x_n \}$ be a sequence in $X \setminus C$ which converges to $x$.

For each $n$ there exists $i_n \geq n$ such that $x_n \in X \setminus C_{i_n}$. Because $C_{i_n}$ is $r$-maximal
maximality implies that $\{x_n\} \cup C_{i_n}$ has diameter larger than $r$ and so there exists $y_n \in C_{i_n}$
for each $n$ such that $d(x_n,y_n) > r$. Let $y$ be a limit point of this sequence in $X_0$. Because $i_n \to \infty$
the point $y$ is in the $limsup = C$.  Furthermore, $d(x,y) \geq r$.
Since $x, y \in C$ we have $d(x,y) \leq r$ and so $x$ and $y$ are antipodal.

If $C$ is $r$-maximal then it is the limit of the constant sequence $C_n = C$ of $r$-maximal sets.

$\Box$ \vspace{.5cm}

\begin{cor}\label{cor1.5}  Let
$C$ be a compact subset of $X$.  Assume that $r > 0$ is a regular value for $C$ (as is always true when the
metric is open).
subset of $X$.

(a) The interior of $C^*_r = \{ d_+(C,\cdot ) \leq r \}$ is $\{ d_+(C,\cdot ) < r \}$.

(b)If the metric is proper and $C$ is an $r$-maximal subset of $X$ then
$A_C$ is a closed, symmetric surjective relation on the compact set $Bdry(C)$.
Furthermore,
\begin{equation}\label{1.6}
C^{\circ} \quad = \quad \bigcap \{  V_r(x) : x \in C \} \quad = \quad \{ d_+(C,\cdot ) < r \}.
\end{equation}
\end{cor}

{\bfseries Proof:} (a): Both $C$ and  its boundary are compact. Let $y \in C$. Because $r$ is a regular value for $y$,
$V_r(y)$ is the interior of $\bar V_r(y)$ and the boundary of each of these sets is
$A_r(y)$. The open set $\{ d_+(C,\cdot ) < r \}$is contained in the interior of $\{ d_+(C,\cdot ) \leq r \} = C^*_r$.
Conversely,  if $x$ is a point of the interior of $C^*_r$ then it is in the interior of
 $\bar V_r(y)$ which is $V_r(y)$. Thus, $x \in V_r(y)$,  for every $y \in C$. By compactness of $C$, $d_+(C,y) < r$.

(b): $A_C$ is always symmetric and
closed.  If $x \in Bdry(C)$ then by Proposition \ref{theo1.4a} there exists $y \in A_C(x)$ and
by Lemma \ref{lem1.1}, $y \not\in C^{\circ}$.  Hence, $A_C$ is a surjective relation on $Bdry(C)$.

Equation (\ref{1.6}) follows from (a) because $C = C^*_r$.

$\Box$ \vspace{.5cm}

Call the metric $d$  \emph{connected} if for every nonempty compact set $C \subset X$
and every $r > 0$ the open set $\{ d_+(C,\cdot ) < r \}$
is connected and if it is non-empty then its closure is $\{ d_+(C,\cdot ) \leq r \}$. In particular,
for any $x \in X$ and $r > 0$ the open set $V_r(x)$ is open and connected with closure $\bar V_r(x)$. Hence, the latter
is a regular closed set which is connected.  Recall that a closed set is called a \emph{regular closed set}
when it is the closure of an open set and so is the closure of its interior. If $r$ is a regular value for $x$ then $V_r(x)$ equals the
interior of $\bar V_r(x)$.

\begin{lem}\label{lem1.6a} Let $X$ have a  connected metric $d$ and let $C$ be a compact subset of $X$.
Assume that $r > 0$ is a regular value for $C$.  If $C_r^*$ has a nonempty
interior then the interior equals $\{ d_+(C, \cdot) < r \}$ and it is a connected open set. $C_r^*$ is a connected,
regular closed set. \end{lem}

{\bfseries Proof:}  By Corollary \ref{cor1.5})(a)  the interior of $C_r^*$
is  $\{ d_+(C,\cdot ) < r \}$ which is connected and nonempty by assumption.  Hence,
$C_r^* = $ \\ $\{ d_+(C,\cdot ) \leq r \}$ is the closure of its interior. Since the closure of a connected set is connected
$C_r^*$ is connected.

$\Box$ \vspace{.5cm}

\begin{theo}\label{theo1.6} Let $X$ be a space with a proper,
 connected metric, let $C$ be a subset of $X$. Assume that $r > 0$ is a regular value for $C$,  as is always true
 when the metric
 is also open.

(a)  If $C$ is $r$-maximal and has a nonempty interior then the interior is connected and $C$ is a connected, regular
closed set.

(b)  If $C$ is $r$-maximal then it is a compact connected set.

(c)  If $C$ is a closed $r$-bounded subset with a nonempty interior
which satisfies the antipodal condition, then $C$ is $r$-maximal.
\end{theo}

{\bfseries Proof:}  (a): Since the metric is proper, $C$ is compact and it equals $C_r^*$ by  Theorem \ref{theo1.3} (f).
 So the result follows
from Lemma \ref{lem1.6a}.

(b):  For any $C$ let $e > 0$ and $r_{e}  = r + 2 e $ and let $C_{e}$ be an $r_{e}$-round
set containing the $r_{e}$-bounded set $ \bar V_{e}(C)$.

Claim:  For every $\epsilon > 0$ there exists a $\delta > 0$ such that $e < \delta$ implies
$C_e \subset V_{\epsilon}(C)$.

Proof:  If there exists $\epsilon > 0$ and a sequence of positive $e_n$'s converging to $0$ and
$x_n \in C_{e_n}$ with $d(x_n,C) > \epsilon$ then by going to a subsequence we obtain a limit point $x$
not in $C$ but such that the diameter of $C \cup \{ x \}$ is $r$. This contradicts $r$-rotundity of $C$.

Thus, we can inductively choose a  sequence of $C_{e_n}$'s with $C_{e_{n+1}} \subset C_{e_n}^{\circ}$ and with
 intersection $C$.  As each $C_{e}$ has
a nonempty interior and is $r_e$ round, it is compact and connected by (a).  Hence, the intersection $C$ is
compact and connected.

(c): Let $D$ be an $r$-maximal set which contains
$C $ and so has a nonempty interior. By (a) $D$ is the closure of its interior and its interior
is connected.   Every point $x$ of the boundary of $C$ has an antipodal point $y$ in $C$. That is $d(x,y) = r$.
Since $r$ is a regular value for $y$, we can apply Lemma \ref{lem1.1}(a) with $C$ replaced by
$D$ to see that $x \in Bdry(D)$. Since the boundary of $D$
is disjoint from the interior of $D$,  the interior of $C$ is nonempty -by hypothesis- and
is clopen in the interior of $D$.  As the latter is connected, the interior of $C$ equals the interior of $D$.
Since $D$ is a regular closed set it must be contained in the closed set $C$ and so $D = C$.  Thus $C$ is
$r$-maximal.

$\Box$ \vspace{.5cm}

I do not know whether an $r$-maximal
set with empty interior can exist when the metric is proper, open and connected. We do have the following  characterization of such odd cases.

\begin{prop}\label{prop1.7}  Let $X$ have a proper, open, connected metric
and let $C$ be an $r$-maximal set with empty interior.
We have $d_+(A_C(x),y) \geq r$ for every $x \in C$ and every $y \in X$. In particular, $A_C \circ A_C (x) = C$ for
every $x \in C$. \end{prop}

{\bfseries Proof:}  Let $U$ be a closed neighborhood of $A_C(x)$ in $C$. Assume that $\{ d_+(U,\cdot ) < r \}$
is nonempty.  Then its closure is $\{ d_+(U,\cdot ) \leq r \}$ which contains all of $C$ and, in particular, contains
$x$.  $d_+(C \setminus U^{\circ},x) < r$ and so by continuity we can choose $y \in X$ with $d_+(U,y) < r$ and
close enough to $x$ that $d_+(C \setminus U^{\circ},y) < r$.  Thus, $d_+(C,y) < r$.  Since $C$ is
 $r$-maximal, (\ref{1.6}) implies that $y$ is in the interior of
$C$ but we have assumed that $C^{\circ} = \emptyset$.

The contradiction implies that for each $y \in X$ and all closed neighborhoods $U$ of $A_C(x)$, $d_+(U,y) \geq r$.
Hence, we can choose a sequence $z_n \in X$ with limit point $z \in A_C(x)$ such that $d(z_n,y) \geq r$.
In the limit we have $d(z,y) \geq r$ and so $d_+(A_C(x),y) \geq r$ for all $y \in X$ and $x \in C$.

In particular, for every $y \in C$, there exists $z \in A_C(x)$ such that $d(y,z) \geq r$ and so $d(y,z) = r$.
That is, $z$ is antipodal to $y$.  Thus, $A_C \circ A_C(x) = C$ for all $x \in C$.

$\Box$ \vspace{1cm}

There is a geometric condition on a metric space which does give us the result we want.

\begin{df} \label{df1.8}Let $X$ be a metric space with metric $d$.

We will say that a point $x \in X$ \emph{lies
between} points $x_0, x_1 \in X$
when for all $y \in X$
\begin{equation}\label{1.7}
d(x,y) \ \leq \ max(d(x_0,y), d(x_1,y)) \hspace{2cm},
\end{equation}
with strict inequality when $d(x_0,y) \not= d(x_1,y) $.

A  \emph{semi-geodesic} between $x_0$ and $x_1$ is a connected subset $G(x_0,x_1)$ of $X$ such
that $x_0, x_1 \in G(x_0,x_1)$ and
 $x $  lies between $x_0$ and $x_1$ for all $x \in G(x_0,x_1) \setminus \{ x_0, x_1 \}$.

We call the metric $d$  \emph{semi-geodesic} when between any two points of $X$ there exists a semi-geodesic.

If the metric is semi-geodesic then we call $C \subset X$  \emph{s-convex}
when any semi-geodesic between two points of $C$
is contained in $C$ and we will call $C$ \emph{w-convex} when between any two points of $C$ some semi-geodesic
between them is contained in $C$.
\end{df}

\begin{prop} \label{prop1.8a} Assume that $X$ has a semi-geodesic metric $d$.

(a) An s-convex set is  w-convex set and a w-convex set is connected.

(b) The intersection of any collection of s-convex sets is s-convex.

(c) If $G(x_0,x_1)$ is a semi-geodesic between distinct points $x_0$ and $x_1$ then $x_0$ and $x_1$ are contained
in the closure of $ G(x_0,x_1) \setminus \{ x_0, x_1 \}$.

(d) For all $x_0, x_1, x \in X$
\begin{equation}\label{1.7a}
x \ \mbox{between} \ x_0 \ \mbox{and} \ x_1 \quad \Longrightarrow \quad max(d(x,x_0),d(x,x_1)) \leq d(x_0,x_1).
\end{equation}
with strict inequality unless $x_0 = x_1$.

Let $x$ lie between $x_0$ and $x_1$. If $x_0 \not= x_1$ then $x \not= x_0$ and $x \not= x_1$.  If $x_0 = x_1$ then
 $x = x_0 = x_1$.

(e) For $x_0, x_1 \in X$ the set
\begin{equation}\label{1.7aa}
\bar G(x_0,x_1) \ =_{def} \ \{ x_0, x_1 \} \cup \{ x \ \mbox{lies between} \ x_0 \ \mbox{and} \ x_1 \}
\end{equation}
is an s-convex, semi-geodesic between $x_0$ and $x_1$. It is  maximal, i.e. it contains every semi-geodesic
between $x_0$ and $x_1$. Furthermore, $\bar G(x_0,x_1) \setminus \{ x_0, x_1 \}$ is s-convex.
\end{prop}

{\bfseries Proof:}  (a): Clearly, s-convex implies w-convex.  If $x_0 \in C$ is fixed and $C$ is w-convex then we can
express $C$ as the union of semi-geodesics $G(x_0,x_1)$ contained in $C$ with $x_1$ varying over $C$.  Hence, $C$ is
connected.

(b): Obvious.

(c):  If $x_0$ were not in the closure of $ G(x_0,x_1) \setminus \{ x_0, x_1 \}$ then $\{ x_0 \}$ would be a
clopen subset of $G(x_0,x_1)$ and it is a proper subset since it does not contain $x_1$. This contradicts
connectedness of $G(x_0,x_1)$.

(d): With $y = x_1$ we have $d(x,x_1) = d(x,y) \leq max(d(x_0,y),d(x_1,y)) = d(x_0,x_1)$.  If  $d(x_0,x_1) = 0$,
then $d(x,x_1) = 0$.  If $d(x_0,x_1) > 0$ then the inequality is strict.  Similarly, for $y = x_0$.

(e): From (d) we can restrict to the case when $x_0 \not=x_1$. Let $z_0, z_1 \in \bar G(x_0,x_1)$. For
 $y \in X$ we have
\begin{equation}\label{1.7ab}
max(d(y,z_0), d(y,z_1)) \quad \leq \quad max( d(y,x_0),d(y,x_1) )
\end{equation}
because each of $z_0, z_1$ is either between $x_0$ and $x_1$ or equal to one of them.

Assume that  $z$ lies between  $z_0$ and $z_1 $ so that
\begin{equation}\label{1.7ac}
d(y,z) \quad \leq \quad max(d(y,z_0), d(y,z_1))
\end{equation}
with strict inequality unless $d(y,z_0) = d(y,z_1)$.  If $z_0 = z_1$ then by (d) again $z = z_0 = z_1$ and so
$z \in G(x_0,x_1)$.  Now assume that $z_0 \not=z_1$.  We show that $z$ lies between $x_0$ and $x_1$.

From (\ref{1.7ab}) and (\ref{1.7ac}) we obtain $d(y,z) \leq max( d(y,x_0),d(y,x_1) )$. If $d(y,z) = max( d(y,x_0),d(y,x_1) )$
then equality holds in (\ref{1.7ab}) and (\ref{1.7ac}) as well.  From equality in (\ref{1.7ac}) we have $d(y,z_0) = d(y,z_1)$
because $z$ lies between $z_0$ and $z_1$. From equality in (\ref{1.7ab}) this common value is $max( d(y,x_0),d(y,x_1) )$.
If either $z_0$ or $z_1$ lies between $x_0$ and $x_1$ then we obtain $d(y,x_0) = d(y,x_1)$. The remaining case is
$\{z_0, z_1 \} = \{x_0, x_1 \}$ so that here $d(y,x_0) = d(y,x_1)$ holds as well because $d(y,z_0) = d(y,z_1)$. As $y$ was arbitrary,
$z$ lies between $x_0$ and $x_1$.

This argument shows that both $\bar G(x_0,x_1)$ and $\bar G(x_0,x_1) \setminus \{ x_0, x_1 \}$ are s-convex and so they are
each connected by (a).
 Hence, $\bar G(x_0,x_1)$ is a semi-geodesic
between $x_0$ and $x_1$.  Because it is s-convex and contains $x_0$ and $x_1$, it contains any semi-geodesic
between them and so is maximal.

$\Box$ \vspace{.5cm}

\begin{lem} \label{lem1.9} Let $X$ be a metric space with a semi-geodesic metric $d$.
\begin{itemize}
\item[(a)] The metric $d$ is connected.
\item[(b)] Let $r > 0$. If $C$ is a nonempty  $r$-bounded subset of $X$ then $C_r^*$ is a closed, s-convex
set. If $C$ is compact then $\{ d_+(C,\cdot) < r \}$ is a nonempty, open, s-convex set with closure $C_r^*$ and
if $C$ is w-convex then $C \cap \{ d_+(C,\cdot) < r \}$ is nonempty.
\end{itemize}
\end{lem}

{\bfseries Proof:}  It is clear from the definitions that the balls $V_r(y)$ and $\bar V_r(y)$ are s-convex.
 $C_r^* = \{ d_+(C, \cdot) \leq r \}  = \bigcap \{ \bar V_r(y) : y \in C \}$ is s-convex and if
 $C$ is compact then $\{ d_+(C,\cdot) < r \} = \bigcap \{  V_r(y) : y \in C \}$ and so it too is s-convex.

(a): Now assume $C$ is compact and that $d_+(C,x_0) < r$ and $d_+(C,x_1) \leq r$.
Let $G(x_0,x_1)$ be a semi-geodesic between $x_0$ and $x_1$.
If $x \in G(x_0,x_1) \setminus \{ x_0, x_1 \}$, and so  $x$ lies between $x_0$ and $x_1$, then for all $y \in C$
 $d(x_0,y) < r$ and $d(x_1,y) \leq r$ and so
 $d(x,y) < r$. Hence, $d_+(C,x) < r$.   By Proposition \ref{prop1.8a}(c)
$x_1$ is in the closure of $ G(x_0,x_1) \setminus \{ x_0, x_1 \}$ and so of $\{ d_+(C,\cdot) < r \}$.
It follows that the metric is connected.

 (b): Now assume that $C$ is compact so that we can write $C$ as the union of compact,
 nonempty subsets $C_1,...,C_N$ each with
 diameter less than $ r $. Choose $x_i \in C_i$ for $i = 1,...,N$. Begin with $y_1 = x_1$.
 Because $C$ is $r$-bounded and the diameter of
 $C_1$ is less than $r$, we have $d_+(C_1,y_1) < r$ and
 $d_+(C,y_1) \leq r$, i.e. $y_1 \in C^*_r \cap C$.

 Assume  we have chosen $y_k \in C_r^*$ with $d_+(C_1 \cup ...\cup C_k, y_k) < r$ for some $k < N$. If
 $y_k = x_{k+1}$ then let $y_{k+1} = y_k = x_{k+1}$. Otherwise, let $G(y_k,x_{k+1})$ be a semi-geodesic between
 the two points.  It is contained in $C_r^*$ by s-convexity. Furthermore, for $z$
 and close enough to $y_k$ is it true that $d_+(C_1 \cup ... \cup C_k, z) < r$ by continuity. On the other hand,
 $d_+(C_{k+1}, x_{k+1}) < r$ and $d_+(C_{k+1}, y_{k}) \leq r$ because $y_k \in C_r^*$.
 By compactness of $C_{k+1}$ it follows that
 $d_+(C_{k+1}, z) < r$ for all $z \in G(y_k,x_{k+1}) \setminus \{ y_k, x_{k+1} \}$. By Proposition \ref{prop1.8a} (c)
 again we can choose
 $y_{k+1} \in  G(y_k,x_{k+1}) \setminus \{ y_k, x_{k+1} \}$  and close enough to
 $y_k$ that
 \begin{equation}\label{1.7b}
 \begin{split}
 d_+(C_1 \cup ... \cup C_{k+1}, y_{k+1}) \qquad =  \hspace{2cm}\\
  max(d_+(C_1 \cup ... \cup C_k, y_{k+1}), d_+(C_{k+1}, y_{k+1})) \ < \ r.
 \end{split}\end{equation}
 Furthermore, if $y_k \in C$ and $C$ is w-convex then we can choose the semi-geodesic $G(y_k,x_{k+1})$
  to lie in $C$ and so
 obtain $y_{k+1} \in C$ as well.

It follows by induction that there exists $x \in X$ such that $d_+(C,x) = d_+(C_1 \cup ... \cup C_N, x) < r$.  Hence,
the open s-convex set $\{ d_+(C, \cdot) < r \}$ is nonempty. Its closure is $C^*_r$ by (a).
Finally, if $C$ is w-convex then the construction yields a point of $C \cap \{ d_+(C,\cdot) < r \}$.

$\Box$ \vspace{.5cm}

\begin{theo}\label{theo1.10}  Let $X$ be a metric space with a proper,  semi-geodesic metric. Let $C$ be a closed,
$r$-bounded subset of $X$ with $r$ a regular value for $C$.  The following are equivalent
\begin{itemize}
\item[i.] $C$ is $r$-maximal.
\item[ii.] $C$ is s-convex and satisfies the antipodal condition.
\item[iii.] $C$ is w-convex and satisfies the antipodal condition.
\item[iv.] $C$ has a nonempty interior and satisfies the antipodal condition.
\end{itemize}

If $C$ is $r$-maximal then it is an s-convex, regular closed subset with interior $\{ d_+(C,\cdot) < r \}$.
\end{theo}

{\bfseries Proof:} (i) $\Rightarrow$ (ii): If $C$ is $r$-maximal then it is
compact because the metric is proper and so by Lemma \ref{lem1.9}(b)
$C = C^*_r$ is  s-convex.
The antipodal condition follows from Theorem \ref{theo1.4a}.

(ii) $\Rightarrow$ (iii): Obvious.

(iii) $\Rightarrow$ (iv): Because $C$ is w-convex, Lemma \ref{lem1.9}(b) implies
there exists a point $x \in C \cap \{ d_+(C,\cdot) < r \}$.
That is, $x \in C$ with no antipode in $C$.  By the antipodal condition every boundary
point of $C$ has an antipode in $C$.  Hence, $x$ is an interior point of $C$.

(iv) $\Rightarrow$ (i): By Lemma \ref{lem1.9}(a) the metric is connected. Hence, if $C$ is $r$-bounded, has a nonempty
interior and satisfies the antipodal condition then it is $r$-maximal by  Theorem \ref{theo1.6}(c).

Applying Lemma \ref{lem1.9}(b) again we see that if $C$ is $r$-maximal then $C = C^*_r$
is the closure of the nonempty open set $\{ d_+(C,\cdot) < r \}$. The latter equals  the
interior by Corollary \ref{cor1.5}(a).

$\Box$ \vspace{1cm}

\section{Minkowski Spaces}

Now we assume that $E$ is a Minkowski space, that is, a vector space with finite
dimension $n$, equipped with a norm $|\cdot|$. Let
$S$ denote the unit sphere in $E$, so that $S = A_1(0)$.  For $x, y \in E$ we will let $[x,y]$ denote
the closed segment $\{ tx + (1-t)y : t \in [0,1] \}$ and similarly we will use $(x,y], [x,y)$ and $(x,y)$
when either or both of the endpoints are omitted.

We first review some elementary facts about such spaces and their convex subsets.

\begin{prop}\label{prop2.1} Let $|\cdot|_0$ be usual Euclidean norm on
${\mathbb R}^n$ with unit sphere $S_0$ and let $|\cdot|$ be
another norm on ${\mathbb R}^n$.  Let $C$ be a nonempty convex subset of ${\mathbb R}^n$.
\begin{itemize}
\item[(a)]  There exist real numbers $M \geq m > 0$ such that for all $x \in {\mathbb R}^n$
\begin{equation}\label{2.1}
m | x |_0  \leq |x| \leq M |x|_0. \hspace{3cm}
\end{equation}
In particular, the topologies and concepts of boundedness induced by the norms $|\cdot|_0$ and $|\cdot|$ agree.

\item[(b)]  If $x_0 \in C$ and $x_1 \in C^{\circ}$ then $(x_0,x_1] \subset C^{\circ}$.

\item[(c)]  If $C$ is closed and bounded and $C^{\circ} \not= \emptyset$ then there is a homeomorphism $h$ on
${\mathbb R}^n$ which maps the unit ball onto $C$.  $C$ is the closure of its
interior and $Bdry(C)$ is homeomorphic to $S$.

\item[(d)]  $S_0$ is homeomorphic to $S$.

\item[(e)]  If $\hat E$ is the affine subspace of $E$ generated by $C$ then the interior
of $C$ with respect to $\hat E$ is
nonempty.

\item[(f)]  The metric on $E$ given by $d(x,y) = |x - y|$ is proper, open, semi-geodesic and connected.
For a subset $C$ of $E$, s-convexity implies convexity implies w-convexity.

\end{itemize}
\end{prop}

{\bfseries Proof:} (a):  If $x = (t_1,...,t_n)$ then $|x| \leq \sum_{i=1}^n \ |t_i| |e_i|$ where $\{ e_1,...,e_n \}$
is the standard basis for ${\mathbb R}^n$. Hence, then norm $|\cdot|$ is a continuous function on ${\mathbb R}^n$ with
respect to the original, Euclidean space topology.  It is positive on the compact sphere $S_0$ and so (\ref{2.1})
holds with $m$ and $M$ the minimum and maximum values on the norm on $S_0$.

(b): By translating if necessary we can assume $x_0 = 0$.  If $V_{\epsilon}(x_1) \in C$ then $V_{t\epsilon}(tx_1) =
tV_{\epsilon}(x_1) \subset C$.

(c): By translating we can assume that $0 \in C^{\circ}$.  Then $r(x) = x/|x|$ defines a continuous map from
$Bdry(C)$ to $S$.  If $u \in S$ then $tu \in C^{\circ}$ for $t$ close to $0$. Since $C$ is closed and bounded,
and hence compact, there exists
$0 < M(u) < \infty$ such that $M(u) = max \{ t : tu \in C \}$. By (b) $tu \in C^{\circ}$ for all $0 \leq t < M(u)$. Hence,
$x = M(u)u$ is the unique point of the boundary of $C$ such that $r(x) = u$. As $r$ is a continuous bijection between
compacta, it is a homeomorphism. Hence, $u \mapsto  M(u) = |r^{-1}(u)|$ is a
continuous, positive real-valued function on $S$.
Let $h(x) =_{def} M(x/|x|)x$ for $x \not= 0$ and $h(0) = 0$ with $h^{-1}(y) = M(y/|y|)^{-1}y$ for $y \not= 0$. Because
$M$ and $M^{-1}$ are bounded, $h$ and its inverse are continuous at $0$. Clearly, $h$ maps the unit ball to $C$.

(d): Since the sphere is the boundary of the unit ball, this follows from (c).

(e): Assume that $\{ x_0,...,x_m \}$ is a maximal affinely independent subset of $C$. Again we translate to assume
$x_0 = 0$ so the $\{ x_1,...,x_m \}$ is a basis for the linear subspace $A$ generated by $C$.  The open simplex
$\{ \sum_{i = 1}^m \ t_i x_i : $ with $ 0 < t_1,...,t_m$ and $ \sum_{i = 1}^m \ t_i < 1  \} $   is a nonempty subset of
$C$ which is open in $\hat E$ (by (a) applied to $\hat E$).

(f): The metric is obviously open.  It is proper because $E$ is finite dimensional. For $x_0, x_1 \in E$ the linear
path $\{ (1-t)x_0 + tx_1 : t \in [0,1] \}$ is a semi-geodesic
because $|x_t - y| \leq (1-t)|x_0 - y| + t|x_1 - y|$ for $t \in [0,1]$.
The metric is connected by  Lemma \ref{lem1.9}(a).  Since a segment between two points is a semi-geodesic it is
clear that the ordinary concept of convexity lies between the concepts of s-convexity and w-convexity.

$\Box$ \vspace{.5cm}

{\bfseries Remark:}  By choosing a basis we can export the results to any normed linear space of dimension $n$.
\vspace{.5cm}

{\bfseries Example:}  Let $E = {\mathbb R}^n$ with $n > 1$ equipped with the $\ell_{\infty}$ norm, $|(a_1,...,a_n)| =
max(|a_1|,...,|a_n|)$.  If $x_0 = (a_1,...,a_n)$ and $x_1 = (b_1,...,b_n)$ with $a_i < b_i$ for $i = 1,...,n$ then
$x = (c_1,...,c_n)$ lies between $x_0$ and $x_1$ iff $a_i < c_i < b_i$ for $i = 1,...,n$. Adjoining $\{x_0, x_1 \} $
to the product of the open intervals $(a_1,b_1),...,(a_n,b_n)$ we obtain the s-convex set $\bar G(x_0,x_1)$ in $E$, see Proposition
\ref{prop1.8a} (e). Hence, the segment $[x_0,x_1]$ is convex but not s-convex. If $x_t$ is a continuous path from $x_0$ to $x_1$ which is
increasing in each coordinate then the image  is a w-convex set but, unless it is the segment, it is not convex.
\vspace{.5cm}

For the rest of this section we assume that $E$ is a Minkowski space.

\begin{theo}\label{theo2.2} Let $C \subset E$ be
an $r$-bounded closed set.
\begin{itemize}
\item[(a)]  $C^*_r = \{ d_+(C,\cdot ) \leq r \}$ is a compact, convex set with a
nonempty interior equal to $\{ d_+(C,\cdot ) < r \}$.

\item[(b)]The following are equivalent
 \begin{itemize}
\item[i.] $C$ is $r$-maximal.
\item[ii.] $C$ is convex and satisfies the antipodal condition.
\item[iii.] $C$ has a nonempty interior and satisfies the antipodal condition.
\end{itemize}

If $C$ is $r$-maximal then it is a convex, regular closed subset with interior $\{ d_+(C,\cdot) < r \}$.
 \end{itemize}
 \end{theo}

 {\bfseries Proof:} Part (a) follows from Lemma \ref{lem1.9} and part (b) follows from Theorem \ref{theo1.10}.

$\Box$ \vspace{.5cm}

Now let $E^*$ denote the dual space of $E$, the space of real linear functionals on $E$.
The space $E^*$ is equipped with the norm $|\omega| = max \{ \omega(x) : x \in S \}$.
We can use the max instead of the
sup because $S$ is compact and the linear functionals are continuous. Thus, for all
$(\omega, x) \in E^* \times E, \ |\omega(x)| \leq |\omega | \cdot |x|$. Let $S^*$ denote the unit sphere in $E^*$.

For any compact convex subset $C$ of $E$ let
\begin{equation}\label{2.2}
G^*_C \quad =_{def} \quad \{ (\omega,x) \in S^* \times Bdry(C) : \omega(x) = max \{ \omega(y) : y \in C \} \},
\end{equation}
which we will call the \emph{Gauss relation}.
For the special case of the unit ball, we define
\begin{equation}\label{2.3}
G^*_1 \quad =_{def} \quad \{ (\omega,x) \in S^* \times S : \omega(x) = 1 \}.
\end{equation}
\vspace{.5cm}

\begin{prop}\label{prop2.3} Let $E$ be a finite dimensional normed linear space with dual space $E^*$.
The relation $G^*_1 \subset S^* \times S$ is a closed, surjective relation.  If $C$ is a compact, convex set
then $G^*_C \subset S^* \times Bdry(C)$ is a closed, surjective relation.  If $C$ is the unit ball in $E$ then
$G^*_C = G^*_1$. \end{prop}

{\bfseries Proof:}  The relation $G^*_1$ is obviously closed. For each $y \in C$ the set $ \{ (\omega,x):
\omega(x) \geq \omega(y) \}$ is obviously closed.  Intersecting over all $y \in C$ we obtain the closed relation
$G^*_C$.

Fix $\omega \in S^*$.  By compactness there exists $x \in C$ at which $\omega$ achieves its
maximum value.  Since a nonzero linear functional is  an open map, this maximum point cannot occur in the interior.
If $C$ is the unit ball then since $|\omega| = 1$ this maximum is $1$ by definition of the norm on $E^*$.
Hence, $G^*_C = G^*_1 $ when $C$ is the unit ball. Since $\omega$ was arbitrary $(G^*_C)^{-1}(S) = S^*$.

Now fix $x \in Bdry(C)$.  There exists $\omega \in S^*$ whose kernel
is the translate of a hyperplane of support through $x$ for the
compact convex set $C$. If $\omega(x)$ is not the maximum on
$C$ then it is the minimum and so $-\omega \in S^*$ takes its maximum at $x$.  Hence, $G^*_C(S^*) = Bdry(C)$.

$\Box$ \vspace{.5cm}

\begin{cor} \label{cor2.4}For $x_1, x_2 \in E, \ |x_1 - x_2| = max \{ \omega(x_1 - x_2) : \omega \in S^* \}.$
\end{cor}

{\bfseries Proof:}  Since $|\omega| = 1, \ \omega(x_1 - x_2) \leq |x_1 - x_2|$ with equality if $x_1 = x_2$. When
the two points are distinct define $x = (x_1 - x_2)/|x_1 - x_2| \ \in S$ and let $\omega \in (G^*_1)^{-1}(x)$ to
get equality.

$\Box$ \vspace{.5cm}

 For any compact subset $C \subset E$ and $\omega \in S^*$ we define
  the \emph{$\omega$ diameter} of $C$ to be
 \begin{equation}\label{2.7}
 \begin{split}
 diam_{\omega}(C) \quad    = \quad   max \{ \ \omega(x_1 - x_2): x_1, x_2 \in C \ \}   \\
 = \quad max \{ \ \omega(x_1): x_1 \in C \ \} \ - \ min \{ \ \omega(x_2):  x_2 \in C \ \} .
\end{split}
 \end{equation}
 From Corollary \ref{cor2.4} it follows that
 \begin{equation}\label{2.8}
 diam(C) \quad = \quad sup_{\omega \in S^*} \ diam_{\omega}(C). \hspace{2cm}
 \end{equation}
\vspace{.5cm}

\begin{theo}\label{theo2.5}Let $C$ be an $r$-bounded subset of  $E$.

(a)For $x_1, x_2 \in C$ and $\omega \in S^*$
\begin{equation}\label{2.4}
\omega(x_1) - \omega(x_2) \quad = \quad r. \hspace{3cm}
\end{equation}
implies that the pair $x_1, x_2$ is antipodal,
$diam_{\omega}(C) = r$ and
\begin{equation}\label{2.5}
\begin{split}
\omega(x_1) \quad = \quad max \{ \omega(x) : x \in C \}, \qquad \mbox{and} \hspace{2cm} \\
\omega(x_2) \quad = \quad min \{ \omega(x) : x \in C \}.  \hspace{3.5cm}
\end{split}
\end{equation}

(b) If  $x_1, x_2 \in C$  is  an antipodal pair, i.e. $|x_1 - x_2| = r$, then
 there exists $\omega \in S^*$ such that $\omega(x_1) - \omega(x_2) = r$.

(c) If  $\omega \in S^*$ with $diam_{\omega}(C) = r$ then there exists an antipodal pair  $x_1, x_2 \in C$
such that $\omega(x_1) - \omega(x_2) = r$.
\end{theo}

{\bfseries Proof:} If $|x_1 - x_2| = r$ then by  Corollary \ref{cor2.4}
 there exists $\omega \in S^*$ such that equation (\ref{2.4}) holds. Conversely, if there exists such an
 $\omega \in S^*$ then by Corollary \ref{cor2.4} $|x_1 - x_2| \geq r$. Equality follows because $C$ is
 $r$-bounded and so $x_1$ and $x_2$ are antipodal.  Clearly, (\ref{2.4}) implies $diam_{\omega}(C) \geq r$ and
 so equality holds because $r \geq diam(C) \geq diam_{\omega}(C)$.

 If $diam_{\omega}(C) = r$ then there exists a pair $x_1,x_2 \in C$ such that (\ref{2.4}) holds.

 Now if equation (\ref{2.4}) holds and $y \in C$ with $\omega(y) \geq \omega(x_1)$ then
 since $C$ is $r$-bounded,
 \begin{equation}\label{2.6}
 r \ \geq \ |y - x_2| \ \geq \ \omega(y) - \omega(x_2) \ \geq \ \omega(x_1) - \omega(x_2) \ = \ r.
 \end{equation}
 Hence, $\omega(y) = \omega(x_1) $ and so $\omega(x_1)$ is the maximum.  Similarly, $\omega(x_2)$ is the minimum.

 $\Box$ \vspace{.5cm}

 \begin{df}\label{def2.6} A subset $C$ has \emph{constant diameter} if
 \begin{equation}\label{2.9}
  diam(C) \quad = \quad  diam_{\omega}(C) \qquad \mbox{for all} \qquad \omega \in S^*.
 \end{equation}
 \end{df}
\vspace{.5cm}

\begin{theo}\label{theo2.7} Assume that a closed, convex set $C$ has  diameter $r$.
\begin{itemize}
\item[(a)] $C$ has constant diameter iff for
every $\omega \in S^*$ there exist  $x_1, x_2 \in C$ such that $\omega(x_1) - \omega(x_2) = r$.
\item[(b)]  If $C$ has constant diameter then $C$ is $r$-maximal.
\end{itemize}
\end{theo}

{\bfseries Proof:} (a): Apply Theorem \ref{theo2.5}.

(b): If $x_1 \in Bdry(C)$ and $\omega \in (G^*_C)^{-1}(x_1)$ then by definition $\omega(x_1)$ is the maximum of $\omega$
on $C$. Let $x_2$ be the point of $C$ where $\omega$ has its minimum.  Then
 $\omega(x_1) - \omega(x_2) = diam_{\omega}(C) = r$. Hence, by Theorem \ref{theo2.5} again $|x_1 - x_2| = r$ and
 so $x_1$ and $x_2$ are antipodal.  As every point of the boundary of $C$ has an antipode,
 $C$ is $r$-maximal by Theorem \ref{theo2.2}.

 $\Box$ \vspace{.5cm}

 For any closed, convex $r$-bounded subset $C$ of $E$ let
\begin{equation}\label{2.10}
\begin{split}
K^*_C \quad =_{def} \quad \{ \  (\omega,x) \in S^* \times Bdry(C) :  \hspace{1cm} \\
 \mbox{for some } \  x_2 \in C \quad
\omega(x) - \omega(x_2) = r \ \}.
\end{split} \end{equation}

\begin{theo} \label{theo2.8} Let $C$ be a closed, convex $r$-bounded subset of $E$.
\begin{enumerate}
\item[(a)]  $K^*_C$ is a closed subset of $G^*_C$.
\item[(b)]  $C$ is $r$-maximal iff $K^*_C(S^*) = Bdry(C)$, ie.  $K^*_C$ projects onto the second coordinate.
\item[(c)] The following conditions are equivalent:
\begin{itemize}
\item[(i)] $(K^*_C)^{-1}(Bdry C) = S^*$.
\item[(ii)] $K^*_C$ is a surjective relation from $S^*$ to $Bdry(C)$.
\item[(iii)] $K^*_C \ = \ G^*_C$.
\item[(iv)] $C$ has constant diameter.
\end{itemize}
\end{enumerate}
\end{theo}

{\bfseries Proof:} (a): It is easy to check that $K^*_C$ is closed. From Theorem \ref{theo2.5} it follows that
$K^*_C$ is contained in $G^*_C$, ie. $\omega(x) - \omega(x_2) = r $ implies $\omega(x)$ is the maximum of $\omega$
on $C$.

(b):  By Theorem \ref{theo2.2}  $C$ is $r$-maximal iff for every $x \in Bdry(C)$ there exists $x_2 \in C$ such that
$|x - x_2| = r$ and by
Corollary \ref{cor2.4} this is true iff  $\omega(x) - \omega(x_2) = r $ for some $\omega \in S^*$ .

(c): Since $G^*_C$ is a surjective relation by  Proposition \ref{prop2.3} it is clear that
(iii) $\Rightarrow $(ii) $\Rightarrow $(i).  Now assume (i) and let $(\omega,x) \in G^*_C$.  By (i) there exist
$x_1, x_2$ such that $\omega(x_1) - \omega(x_2) = r$. By (a) $(\omega, x_1) \in G^*_C$ and so the maximum of
$\omega$ on $C$ is  $\omega(x_1) = \omega(x)$.  Hence, $\omega(x) - \omega(x_2) = r $ and so $(\omega,x) \in K^*_C$.
Thus, (i), (ii) and (iii) are equivalent.

By Theorem \ref{theo2.7}(a) (iv) is equivalent to (i).

$\Box$ \vspace{.5cm}

{\bfseries Remark:}  Notice that if $\omega(x_1) - \omega(x_2) = r$ then $(\omega,x_1) \in K^*_C$ and
$(- \omega,x_2) \in K^*_C$.
\vspace{.5cm}

{\bfseries Example:}  When $C$ is the unit ball, $G^*_C = G^*_1$ by  Proposition \ref{prop2.3}. It follows that
for every $\omega \in S^*$ the maximum value on $S$ is $1$ and the minimum value is $-1$. Hence, the unit ball
has constant diameter $2$.  It follows that any ball in $E$ of radius $r/2$ has constant diameter $r$ and is
$r$-maximal.

If $n = 1$ then $E$ is isometric to ${\mathbb R}$. In a one-dimensional
normed linear space  the closed balls of radius $r/2$, i.e. the closed intervals of length $r$, are all of
the $r$-maximal sets.

For any $n$ if $E$ is ${\mathbb R}^n$ equipped with the $\ell_{\infty}$ norm: $|(x_1,...,x_n)|
= max(|x_1|,...,|x_n|)$ then the closed balls of radius $r/2$ are all of the $r$-maximal sets. To see this,
project the $r$-maximal set $C$ to coordinate $i$.  The image is contained in some interval $I_i$ of length $r$.
The product of the intervals $I_1 \times ... \times I_n$ is a ball of radius $r/2$ and so is $r$-maximal. The
product contains $C$ and so, by maximality, equals $C$.
\vspace{.5cm}

 $C$ is called \emph{strictly convex}
if whenever $x_0, x_1$ are two distinct elements of $C $ then $(x_0,x_1) \subset C^{\circ}$. A nonempty, strictly convex
set is either a singleton or else has a nonempty interior.

\begin{prop} \label{prop2.9} If $C$ is a strictly convex, closed, bounded set then  every $\omega \in S^*$
achieves its maximum value on $C$ at a unique point which lies in $Bdry(C)$. The relation $G^*_C$ is a continuous
function from $S^*$ onto $Bdry(C)$. \end{prop}

{\bfseries Proof:}  If $a = \omega(x_0) = \omega(x_1)$ with $x_0 \not= x_1$ in $C$ then $a = \omega(x_t)$
with $x_t = tx_1 + (1-t)x_0$ for all $t \in [0,1]$.  By strict convexity, $x_t \in C^{\circ}$ for $t \in (0,1)$.
Because $\omega$ is an open map to ${\mathbb R}$ the value $a$ cannot be a maximum.

Hence, the closed relation $G^*_C$ is a function from $S^*$ to $Bdry(C)$.  It is easy to check that a function
which is closed relation between compacta is continuous, see Akin (1992) Corollary 1.2.

$\Box$ \vspace{.5cm}

\begin{cor} \label{cor2.10} If $C$ is a strictly convex, closed, $r$-bounded subset then $K^*_C$ is a continuous
function on a closed subset $dom(K^*_C) \subset S^*$. If $\omega \in dom(K^*_C)$ then $- \omega \in dom(K^*_C)$.

$C$ is $r$-maximal iff $K^*_C$ is a surjective map from $dom(K^*_C)$ to
$Bdry(C)$ and $C$ has constant diameter iff $dom(K^*_C) = S^*$. \end{cor}

{\bfseries Proof:}  The statements  in the first paragraph are immediate from
Proposition \ref{prop2.9} and the Remark following
Theorem \ref{theo2.8}.  The domain of the
function $K^*_C$ is the closed subset
$dom(K^*_C) = $ \\ $(K^*_C)^{-1}(Bdry(C))$.

 The second paragraph follows from Theorem \ref{theo2.8} itself.

$\Box$ \vspace{.5cm}

The norm $|\cdot|$ on $E$ is called \emph{strictly convex} when the unit ball $V_1(0)$ is strictly convex.
Since we obtain any other ball by dilating and translating it then follows that every closed ball is strictly
convex.

\begin{theo} \label{theo2.11} Assume that $E$ has a strictly convex norm.
If $C$ is a compact subset of $E$ and $r > 0$ then $C^*_r = \{ d_+(C,\cdot) \leq r \}$ is strictly convex.
In particular, if $C$ is
$r$-maximal then it is strictly convex. \end{theo}

{\bfseries Proof:}  If for distinct points $x_1, x_2 \ d_+(C,x_1), d_+(C,x_2) \leq r$
then for every $y \in C$ we have $x_t \in V_r(y)$ for each $x_t \in (x_1,x_2)$. By compactness
$d(y,x_t) < r$ for all $y \in C$ implies $d_+(C,x_t) < r$ and so $x_t \in \{ d_+(C,\cdot) < r \}$ which is
the interior of $C^*_r$.

$\Box$ \vspace{.5cm}

We conclude this section with some results we will use later.

\begin{prop} \label{prop2.12} Assume that $C$ is a bounded closed set in $E$ with a nonempty interior.
If for every point $x \in Bdry(C)$
there exists a hyperplane of support for $C$ containing $x$ then $C$ is convex. \end{prop}

{\bfseries Proof:}  Let $U = C^{\circ}$ so that $U$ is nonempty. The assumption
about hyperplanes of support says that for every
$x \in Bdry(C)$ there exists $\omega_x \in S^*$ which achieves its maximum on $C$ at $x$.  Let
\begin{equation}\label{2.11}
\hat C \quad =_{def} \quad \bigcap_{x \in Bdry(C)} \  \{ y  \in E : \omega_x(y) \leq \omega_x(x) \ \}.
\end{equation}

$\hat C$ is a closed, convex set which contains $C$ and so has nonempty interior.  Recall that each $\omega \in S^*$
is an open map. Hence, for each
$x \in Bdry(C)$ there exist $z \in E$ arbitrarily close to $x$ with $\omega_x(z) > \omega_x(x)$ and so with
$z \not\in \hat C$.  Thus, $Bdry(C) \subset Bdry(\hat C)$. Since the boundary of $\hat C$ is disjoint from its
interior it follows that $U$ is a clopen subset of the interior of $\hat C$.  By Proposition \ref{prop2.1} (c)
the interior of $\hat C$ is homeomorphic to the open unit ball and so is connected.  Hence, $U$ is the interior of
$\hat C$. Because $\hat{C}$ is the closure of its interior, we have
\begin{equation}\label{2.12}
\hat{C} \ = \ \overline{U} \ \subset \ C \ = \ U \cup Bdry(C) \ \subset \ \hat{C}.
\end{equation}
Thus, $C = \hat{C}$.

$\Box$ \vspace{1cm}

\begin{theo}\label{theo2.13} Let $\{C_n \}$ be a sequence of compact subsets of $E$ which converge to
$C$ in the Hausdorff metric.  If each $C_n$ is $r$-maximal then $C$ is $r$-maximal.  If each $C_n$ is a convex
set with constant diameter $r$ then $C$ is a convex set with constant diameter $r$. \end{theo}

{\bfseries Proof:}  Just as with (\ref{1.1e}) it is easy to check that for any $\omega \in S^*$ and compacta $C,D$:
\begin{equation}\label{2.13}
\begin{split}
|diam(C) - diam(D)| \quad  \leq \quad 2 d(C,D), \\
|diam_{\omega}(C) - diam_{\omega}(D)| \quad  \leq \quad 2 d(C,D)
\end{split}
\end{equation}
Hence, the conditions $diam(C) \leq r$, $diam(C) = r$ and $diam_{\omega}(C) = r$ are preserved by limits in the
Hausdorff metric.

The function $q : E \times E \times [0,1] \to E$ by $(x_0,x_1,t) \mapsto (1-t)x_0 + tx_1$ is continuous and
so the association $C \mapsto q(C \times C \times [0,1])$ is continuous with respect to the Hausdorff metric
(see, e.g., Akin (1993) Proposition 7.16).
Hence, the set of $\{ C : C = q(C \times C \times [0,1]) \}$ is a closed set. This is the collection of convex sets.

It follows that if each $C_n$ is a convex set of constant diameter $r$ then so is $C$.

Now if each $C_n$ is $r$-maximal, then each is convex and so $C$ is convex.  By  Theorem \ref{theo1.4a} $C$
satisfies the antipodal condition. By Theorem \ref{theo2.2} (c) $C$ is $r$-maximal.

$\Box$ \vspace{.5cm}

It is a classic result, which we will reprove in the next section, that for a Euclidean an $r$-bounded subset is
$r$-maximal iff it has constant width.  That is, in the Euclidean case the converse of Theorem \ref{theo2.7} (b)
holds. Eggleston (1965)  gave examples of Minkowski spaces with $r$-maximal subsets which are not of constant width.
\vspace{1cm}

\section{Euclidean Spaces}

From now on we assume that $E$ is a Euclidean space of dimension $n \geq 2$.  That is, $E$ is an $n$ dimensional
linear space with norm given by an inner product an so isometric to ${\mathbb R}^n$ with the usual metric. If
$x_1, x_2$ are distinct points of norm $1$ then
\begin{equation}\label{3.1}
\begin{split}
| tx_1 + (1-t)x_2 |^2 \ = \ t^2|x_1|^2 + 2t(1-t) x_1 \cdot x_2 + (1-t)^2|x_2|^2 \\
= \quad 1 - t(1-t)(1 - \cos \theta ) \hspace{4cm}
\end{split}
\end{equation}
which is less than $1$ since the angle $\theta$ between $x_1$ and $x_2$ is not $0$.  That is, the norm is
strictly convex. It follows from Theorem \ref{theo2.11} that for $C \subset E$ compact and $r > 0$ the
dual $C^*_r$ is strictly convex.  However, in this case we can do better.

Let $y_1$ and $y_2$ be two distinct points of the Euclidean space $E$.
For $x \in X$ with $d(x,y_1) = d(x,y_2) = r$
there is a unique circle in $X$ which is centered at $x$ and which passes through $y_1$ and $y_2$.
The  \emph{arc $\a$ between $y_1$ and $y_2$ with center $x$}  is the smaller of the two arcs
of this circle with endpoints $y_1$ and $y_2$. If $x$ is the midpoint of $[y_1,y_2]$ then either
semicircle with endpoints $y_1$ and $y_2$ is considered to be ``the" arc $\a$. The arc $\a$ is said to have
radius the common distance $r$. We call $S \setminus \{y_1, y_2 \}$ the
\emph{open arc between $y_1$ and $y_2$ with center $x$}.

\begin{df}\label{def3.1} Let $E$ be a Euclidean space  of dimension at least 2 and let $r > 0$. A subset $C $ of $X$ is
called \emph{$r$-convex} if it is convex and $\a \subset C$ whenever $\a$ is an arc of radius at
least $r$ between  two points
of $C$.\end{df}
\vspace{.5cm}

Notice that as $r$ tends to infinity the arcs between $y_1$ and $y_2$ tend to the segment $[y_1,y_2]$.
From this it is not hard to show that convexity actually follows from the arc condition.
Clearly the closure of an $r$-convex set is
$r$-convex and any intersection of $r$-convex sets is $r$-convex.  For $C \subset E$
 the $r$-convex hull of $C$ is the intersection of all $r$-convex sets which contain $C$ or, equivalently, the
 minimum $r$-convex set containing $C$.

\begin{prop}\label{prop3.2}  If $C$ is  $r$-convex  then it is strictly convex. \end{prop}

{\bfseries Proof:} For $y_1, y_2$ two distinct points of $C $ with $d(y_1,y_2) = 2 \ep $ let $y$ be the
midpoint of the segment between them and $P$ be the hyperplane through $y$ perpendicular to the line through
them. Choose $a > max(r, \ep)$.  Let $S$ be the sphere in $P$ centered at $y$
with radius $a - \sqrt{a^2 - \ep^2}$. For each point  of $  S$
lies on an arc between $y_1$ and $y_2$ of radius $a > r$.  So by $r$-convexity $S \subset C$.  By convexity
$C$ contains the convex hull of $S \cup \{ y_1, y_2 \}$ whose interior contains every point on the open
segment between $y_1$ and $y_2$.  Thus, $C$ is strictly convex.

$\Box$ \vspace{.5cm}

We now require a bit of plane geometry.

\begin{lem}\label{lem3.3} Assume $P$ is a Euclidean plane. Let $\a$ be an arc
between distinct points $y_1, y_2 \in X$ with center $x$ and radius $r$. Let $B$ be the
ball of radius $r$ centered at $x$ so that $\a \subset S,$ the boundary circle of $B$.
Let $B_0$ be a ball of radius
$r_0 > 0$ and boundary circle $S_0$. If $\a \subset B_0$ and $\a \cap S_0$ contains a point of the open arc
$ \a \setminus \{y_1, y_2 \}$ then $r_0 \geq r$ and if $r_0 = r$ then $x$ is the center of $B_0$
and so $S_0 = S \supset \a$.
\end{lem}

{\bfseries Proof:} Let $y \in S_0 \cap \a \setminus \{y_1, y_2 \}$. If $\a$ intersects $S_0$
transversely at $y$  then
part of $\a$ near $y$ on one side or the other lies outside of $B_0$ contra assumption. Hence, $\a$
intersects $S_0$ tangentially at $y$. The center $x$ of $\a$ lies on the line through $y$  of the
 perpendicular to the tangent of $\a$ at $y$. The center $x_0$ of
$S_0$ must lie on the same line. If $x$ and $x_0$ lay on opposite sides of $y$ on the line,
then $B \cap B_0 = \{ y \}$ again contradicting $\a \subset B_0$. If
$x_0$ is on the open segment between $x$ and $y$ (i.e. $r_0 < r$) then $B_0 \subset B$ with
$B_0 \cap S = \{ y \}$ and so $B_0 \cap \a = \{ y \}$ which is still false. Hence, $r_0 \geq r$.
Furthermore, if $r_0 = r$ then $x_0 = x$, $B_0 = B$ and $S_0 = S$.

$\Box$ \vspace{.5cm}

\begin{lem}\label{lem3.4} A closed ball $B$ in $E$ of radius $r >0 $  is $r$-convex.
\end{lem}

{\bfseries Proof:} Since $B$ closed and the arcs of radius $a$ move continuously with the endpoints it
suffices to consider the case when $y_1, y_2$ are in the interior of $B$. For any arc between $y_1$ and $y_2$
the three points $y_1, y_2$ and the center $x$ determine a plane $P$ in $E$. Any plane  through $y_1$ and
$y_2$ intersects $B$ in a disc $B_1$ with radius $r_1 \leq r$. Let $x$ move along the line $L$ in
$P$ which is the perpendicular bisector of the segment $[y_1,y_2]$ and for each such $x$ let $\a$ be the arc
between $y_1$ and $y_2$ centered at $x$.

As $x$ moves toward infinity the arc $\a$ approaches the segment $[y_1,y_2]$ which is contained in the
 interior of $B_1$ with respect to the plane $P$. As $x$
 moves from infinity toward the mid-point either all of the arcs remain in the
interior of $B_1$ or there is a first position $x^*$  such that the arc touches the boundary circle of $B_1$.
This arc in $P$ centered at $x^*$ is contained in $B_1$ and by Lemma \ref{lem3.3}
has radius less than $r_1$ and so less than $r$.

Hence, as long as $d(x,y_1) = d(x,y_2) \geq r$ the arc is contained in the interior of $B_1$.

$\Box$ \vspace{.5cm}

From the proof we obtain:

\begin{prop} \label{prop3.4a}  For two distinct points $x_0$ and $x_1$ in a Euclidean space $E$ the only
points of $E$ which lie between $x_0$ and $x_1$ in the metric space sense are those of the open
segment $(x_0,x_1)$ and so  $ [x_0,x_1]$ is the unique semi-geodesic between $x_0$ and $x_1$. In particular, for a
Euclidean space the concepts of s-convexity, convexity and w-convexity agree.
\end{prop}

{\bfseries Proof:}  To show that no point $x$ of $E \setminus (x_0,x_1)$ lies between $x_0$ and $x_1$ it
suffices to restrict to the plane through $x_0, x_1$ and $x$ and so to reduce to dimension two.

Let $L$ be the line through $x_0$ and $x_1$ and let $\bar L$ be the perpendicular bisector of $[x_0,x_1]$.
For any $y \in \bar L$ the in-between point $x$ must lie in the disc centered at $\bar V_r(y)$ with
$r = |x_0 - y| = |x_1 - y|$. As $y$ moves out toward infinity in one direction, the intersections of the
discs with the closed half-space on the other side decrease with intersection $[x_0,x_1]$. So if $x$ is
in the closed half-space then it must lie in $[x_0,x_1]$.

It follows that $[x_0,x_1]$ is the maximal semi-geodesic $\bar G(x_0,x_1)$.  On the other hand, removing
any points other than $x_0$ and $x_1$ disconnects the interval and so no proper subset of $[x_0,x_1]$ is
a semi-geodesic between $x_0$ and $x_1$.

The convexity results are then obvious.

$\Box$ \vspace{.5cm}

\begin{prop} \label{prop3.5} For $C $ a bounded subset of $E$, let
$D$ be the $r$-convex hull of $C$.  The set $C^*_r$ is $r$-convex and $C^*_r = D^*_r$. \end{prop}

{\bfseries Proof:} By Lemma \ref{lem3.4} the ball $\bar V_r(y)$ is $r$-convex and so it contains $D$ iff it contains $C$.
Hence, $C^*_r = D^*_r$. Furthermore, $C^*_r$ is the intersection of $\{ \bar V_r(x) : x \in C \}$ an so is
$r$ convex.

$\Box$ \vspace{.5cm}

\begin{theo}\label{theo3.6} If $C$ is an $r$-maximal set  in $E$, then $C$ is $r$-convex.
 \end{theo}

 {\bfseries Proof:}  By  Theorem \ref{theo1.3}(f) $C = C^*_r$. So $C$ is $r$-convex by Proposition \ref{prop3.5}

 $\Box$ \vspace{.5cm}

\begin{cor}\label{cor3.7} Let $C$ be an $r$-maximal set in $E$.
Assume that for $x \in C$ there are two distinct points $ y_1, y_2 \in C$ antipodal to $x$. Let $\a$ the  arc  between
$y_1$ and $y_2$ with center $x$. The angle subtended by $\a$ is at most $\pi/3$ and every point of $\a$ is a point of
$Bdry (C)$ antipodal to $x$.
\end{cor}

{\bfseries Proof:} The distance $d(y_1,y_2) \leq r$ which equals the radius of the arc.  Since the chord
between the end-points has length at most the radius, the angle subtended is at most $\pi/3$.

By Theorem \ref{theo3.6} $\a \subset C$. Since $x$ is the center of $\a$, $d(y,x) = r$ for every $y \in \a$
and so every point of $\a$ is antipodal to $x$. By Lemma \ref{lem1.1} $\a \subset Bdry (C)$.

$\Box$ \vspace{.5cm}

\begin{theo}\label{theo3.8} Let $C$ be an $r$-maximal subset in $E$.
Assume that $ y_1, y_2 $ are distinct points of $E$ and that $\a$ is an  arc  between them
with radius $r$ and center $x$.  If $\a \subset Bdry(C)$ then $x \in C$ and for each $y$ in
the open arc $\a \setminus \{ y_1, y_2 \}, A_C(y) = \{ x \}$.  That is, $x$ is the unique
point of $C$ antipodal to $y$. \end{theo}

{\bfseries Proof:}  Let $z \in Bdry(C)$ be a point antipodal to $y$. We will show that $z = x$.

Let
$P$ be the  plane which contains $y_1, y_2$ and $x$ and let $\hat E$ be smallest affine
subspace which  contains $P$ and $z$. Let $B$ be the ball in $\hat E$ with center $z$ and radius $r$.
Because $\a \cup \{z \} \subset C$ and $diam(C) = r$ we have $\a \subset B$. Let $S$ be the boundary
sphere of $B$ in $\hat E$.

First we show that $z \in P$ and so $\hat E = P$. If not then the intersection $P \cap B$ is a disk in
the plane $P$ with radius less than $r$  and with boundary circle $P \cap S$.
Furthermore, $\a \subset P \cap B$ and $\a \cap (P \cap S)$ contains the point $y$ of the open arc.
By Lemma \ref{lem3.3} this can't happen.

Hence, $z \in P$ and so $B$ is a disk in $P$ with radius $r$ which contains $\a$ and whose boundary
circle meets the point $y$ of the open arc.  By Lemma \ref{lem3.3} again the center $z$ of $B$ agrees
with the center $x$ of $\a$.  That is, $z = x$ as required.

$\Box$ \vspace{.5cm}

{\bfseries Example:}  With $P$ the Euclidean Plane let $C_0$ be a regular polygon with $2k+1$ vertices
($k \geq 1$) indexed, in order, by the additive group ${\mathbb Z}/(2k+1){\mathbb Z}$. Let $r =
d(v_i,v_{i+k}) = d(v_i,v_{i-k})$. For other pairs of vertices, $d(v_i,v_j) < r$. If $C_1$ is an
$r$-maximal set containing these vertices (and hence by convexity the polygon $C_0$ and its interior)
then $C_1$ contains the arc between $v_{i+n}$ and $v_{i-n}$ centered at $v_i$ which consists of points
antipodal to $v_i$. The convex hull $C$ of these arcs is a closed subset such that every point of the boundary
has antipodal points. Hence, by Theorem \ref{theo2.2}(c) $C$ is $r$-maximal and contains the vertices.  By maximality
$C = C_1$.  Thus, $C$ is the unique $r$-maximal set which contains the vertices. These are the Relleaux polyhedrons.
 With $k = 1$ it is the classical
Relleaux triangle.

The vertices need not be evenly spaced.  Beginning with the regular polygon with $k \geq 2$ define the
$2k+1$ pointed star to be the graph with vertex $v_i $ connected to $v_{i+k}$ and to $v_{i-k}$ by segments
of length $r$. The star is not rigid and so we can vary the vertices slightly and still have a set of
diameter $r$. Connecting $v_{i+k}$ and $v_{i-k}$ in the new positions by the arc centered at the new position
for $v_i$ we get the unique $r$-maximal set containing the perturbed vertices.
\vspace{.5cm}

Now we show that in the Euclidean case the $r$-maximal sets are exactly the sets of constant diameter $r$, see, e.g.
Bonnesen and Fenchel (1987) and Eggleston (1958). For
any Minkowski space Theorem  \ref{theo2.7} says that a set of constant diameter $r$ is
$r$-maximal.  We prove the converse in the Euclidean case.

In the Euclidean case we can replace the use of the dual space $E^*$ by using $E$ instead. The inner product
provides  the Riesz Representation isometry from $E $ to $E^*$ by $v \mapsto \omega_v$ where
$\omega_v(x) =_{def} v \cdot x$. On the unit spheres the isomorphism restricts to:
\begin{equation}\label{3.2}
G^*_1 \quad = \quad \{ (\omega_u,u) : u \in S \}
\end{equation}
because for unit vectors $u, u_1$ we have $u \cdot u_1 = 1 $ iff  $u_1 = u$. Using this we define for  $C$ a
compact convex set the relation $G_C$ dual to $G^*_C$:
\begin{equation}\label{3.3}
G_C \quad = \quad \{ (u,x) \in S \times C : u \cdot x = max \{ u \cdot y : y \in C \} \ \}.
\end{equation}

\begin{prop} \label{prop3.10}For a compact convex set $C$ the relation
$G_C$ is a closed surjective relation from $S$ to $Bdry(C)$.
If $C$ is strictly convex then $G_C$ is a continuous surjective function from $S$ to $Bdry(C)$. \end{prop}

{\bfseries Proof:}  Since the association $u \mapsto \omega_u $ is a
homeomorphism from $S$ to $S^*$ the first result follows
from  Proposition \ref{prop2.3} and the second from Proposition \ref{prop2.9}.

$\Box$ \vspace{.5cm}

Thus, when $C$ is a compact strictly convex set, e.g. when $C$ is $r$-maximal, we can write $x = G_C(u)$ for the
unique  $x \in Bdry(C)$ such that $(u,x) \in G_C$. That is, given $u \in S, \ x = G_C(u)$ is the unique point of
$C$ at which $z \mapsto \omega_u(z) = u \cdot z$ takes its maximum on $C$.  Of course, $G_C(u) \in Bdry(C)$.
Following Lachand-Robert and Oudet we will call this function
the \emph{inverse Gauss map}.

Now we need another bit of plane geometry.

\begin{lem}\label{lem3.11} Let $x_1, x_2$ be distinct points of a Euclidean plane $P$ with $|x_1 - x_2| = r$.
Let $L_1, L_2$ be the
lines through $x_1$ and $x_2$ respectively which are perpendicular to  the segment $[x_1,x_2]$ and so are parallel.
Let $B_1, B_2$ be the discs or radius $r$ centered at $x_1$ and $x_2$ respectively, with boundary circles
$S_1, S_2$. Let $U$ be the open strip consisting of the points of $P$ between the lines $L_1$ and $L_2$.

Assume that $y \in B_1 \cap U$ and that $\a$ is the arc between $y$ and $x_1$ of radius $r$ with center on $S_1 \cap U$.
If $\a \subset \overline{U}$ then $y \in B_2$.
\end{lem}

{\bfseries Proof:}  Choose coordinates so that $x_2$ is the origin and $x_1 = (0,r)$. Assume that $y = (a,b)$ with
$0 < b < r$.  If $a = 0$ then $y \in [x_1,x_2] \subset B_2$.  Now assume that $y \not\in B_2$.

Without loss of generality we can assume that $a > 0$.    Moving along the segment $[y,x_1]$ from $y$ to $x_1$
let $x_3$ be the first entrance into $B_2$.  Thus, $[x_3,x_1]$ is a chord of the circle $S_2$ and its perpendicular
bisector $L_3$ passes through the center $x_2$. The perpendicular bisector of $ [y,x_1]$ is parallel to $L_3$ but
to the right of it and so it intersects the semi-circle $S_1 \cap U$ at a point $(c,d)$ with $c > 0$.  Then
the arc $\a$ which connects $y$ and $x_1$ with center $x_4 = (c,d)$ is tangent at $x_1$ to the line perpendicular
to $[x_4,x_1]$ and is above the tangent line $L_1$ to $S_1$ near
 $x_1$. In particular, $\a$ is not contained in $\overline{U}$. Contrapositively,
 $\a \subset \overline{U}$ implies $y \in B_2$.

$\Box$ \vspace{.5cm}

\begin{theo} \label{theo3.12} Assume that $C$ is an $r$-bounded subset of $E$,

(a) If $x_1, x_2$ is an antipodal pair in $C$ i.e.
$|x_1 - x_2| = r$, and  $u = (x_1 - x_2)/r \in S$ then $\omega_u \in S^*$ with $\omega_u(x_1) - \omega_u(x_2) = r$
and $\omega_u(x_1) = max \{ \omega_u(y) : y \in C \}$.
The hyperplane through $x$  which is perpendicular $u$
 is a hyperplane of support for $C$ which contains $x$.  In particular, if $x$ has more than one
antipodal point in $C$ then $C$ has more than one hyperplane of support through $x$. If $\ C$ is convex then
$(u,x_1) \in G_C$.

(b) Assume that $C$ is an $r$-maximal subset of $E$.  If $u \in S$ and $x_1 = G_C(u)$ then
$x_2 = x_1 - r u = G_C(-u)$. The pair $\{ x_1, x_2 \}$ is the unique antipodal pair in $C$ such that the
segment $[x_1,x_2]$ is parallel to $u$.
\end{theo}

{\bfseries Proof:} (a): Direct computation shows $\omega_u(x_1) - \omega_u(x_2) = r$. The rest follows from
Theorem \ref{theo2.5} and the definition of $G_C$ in the convex case.

(b): This is much more delicate.  The main issue is to  show that $x_2 \in C$.  It will suffice to
show that $C \subset \bar V_r(x_2)$ and so $x_2 \in C^*_r$, because $C = C^*_r$ by Theorem
\ref{theo1.3} (f).

 Now let $y \in C$. If $y = x_1$ or $x_2$
then $y \in \bar V_r(x_2)$.  Now assume that $y$ is distinct from $x_1$ and $x_2$ and
let $P$ be the plane in $E$ which contains $x_1, x_2, y$. We apply the lemma.

Let $P \cap V_r(x_i) $ be the disc $B_i$ centered at $x_i$ with radius $r$ for $i = 1,2$. The strip $U =
\{ w \in P : \omega_u(x_2) < \omega_u(w) < \omega_u(x_1) \ \}$ . Because $x_1 \in C$ and $C$ is $r$-bounded,
$C \subset V_r(x_1)$ and so $P \cap C \subset B_1$.  Also for $w \in C \ $ $ \omega_u(w) \leq \omega_u(x_1)$ with
equality only with $w = x_1$ by strict convexity of $C$.  Hence, $P \cap C = \overline{U} \cap C \subset
\{ x_1, x_2 \} \cup U \cap C$.  Hence, $y \in U \cap C \subset U \cap B_1$.  Let $\a$ be the arc of radius $r$
which connects $x_1$ and $y$ with center on a point of $Q_1 \cap U$. Because $C$ is $r$-convex
it follows that the arc $\a \subset C \cap P \subset \overline{U}$.  Then Lemma \ref{lem3.11} implies that
$y \in B_2$ and so $y \in \bar V_r(x_2)$ as required.

Clearly, $r = |x_1 - x_2| = \omega_u(x_1) - \omega_u(x_2)$ so that $\{ x_1, x_2 \} $ is an antipodal pair with
$[x_1,x_2]$ parallel to $u$. Furthermore, by Theorem \ref{theo2.5} $u \cdot x_2 = min \{ u \cdot y : y \in C \}$ and
so $x_2 = G_C(-u)$. If $[y_1,y_2]$ is parallel to $u$ then $(y_1 - y_2)/r = \pm u$ and by renumbering if necessary
we can assume that the sign is positive so that $y_2 = y_1 - ru$. Since $\omega_u(y_1) - \omega_u(y_2) = r$,
Theorem \ref{theo2.5} again implies $\omega_u(y_1)$ is the maximum value of $\omega_u$ on $C$.  But $x_1$ is the
unique point $G_C(u)$ where $\omega_u$ takes its maximum.  Thus, $y_1 = x_1$ and $y_2 = x_2$.

$\Box$ \vspace{.5cm}

\begin{cor}\label{cor3.13}  If $C$ is an $r$-maximal subset of the Euclidean space $E$
then $C$ has constant diameter $r$. Furthermore, $C$ is the union of the antipodal segments in $C$.
That is, if $ \ x \in C$ then there exists an antipodal pair  $\ x_1, x_2$ in $Bdry(C)$ such that $\ x \in [x_1,x_2]$.
\end{cor}

{\bfseries Proof:}  Given $\omega \in S^*$
let $u $ be the unit vector of $E$ such that $\omega = \omega_u$. With $x_1 = G_C(u)$ and $x_2 = x_1 - ru$
Theorem \ref{theo3.12} (b) says that $x_1,x_2$ is an antipodal pair in $C$ with $\omega(x_1) - \omega(x_2) = r$.
Since $\omega \in S^*$ was arbitrary
Theorem \ref{theo2.7}(a)  implies that $C$ has constant diameter $r$.

If $x \in C$, let $r_0 = d_+(C,x) $ so that $r \geq r_0 > 0$.  Let $x_1 \in C$ with $|x - x_1| = r_0$.
Thus, $C$ is contained in the ball $V_{r_0}(x)$ and $x_1$ is on its boundary sphere. Let $u = (x_1 - x)/r_0$.
$\omega_u$ achieves its maximum on $V_{r_0}(x)$ at $x_1$ and so achieves its maximum on $C$ at $x_1$.  That is,
$G_C(u) = x_1$.  By Theorem \ref{theo3.12} (b) again $x_2 = x_1 - ru$ is an antipodal point for $x_1$.
Since $r_0 \leq r$,
$x \in [x_1,x_2]$.

$\Box$ \vspace{1cm}

\section{Parametrizations}

Now let $E$ be the Euclidean space ${\mathbb R}^n$ with $n \geq 2$ and equipped with
the usual metric.

First assume that $C$ is a $r$-maximal subset of $E$, or, equivalently, a closed, convex set with constant diameter $r$.

In the previous section we defined the inverse Gauss map  $G_C : S \to Bdry(C)$ which associates to each
$u \in S$ the unique point of $C$ at which $\omega_u$ achieves its maximum. By Proposition  \ref{prop3.10} this is
a continuous, surjective function because an $r$-maximal set is strictly convex. By
Theorem \ref{theo3.12}(b), if $x_1 = G_C(u)$ then $x_2 = x_1 - ru = G_C(-u)$ is the point of $C$  antipodal to
$x_1$ with $[x_1,x_2]$ parallel to $u$.  Thus we have
\begin{equation}\label{4.1}
u \quad = \quad (G_C(u) \ - \ G_C(-u))/r.  \hspace{3cm}
\end{equation}
By Theorem \ref{theo3.12} (a) as $u$ varies over $S$, the pairs $G_C(u),G_C(-u)$ vary over all antipodal pairs.
from Corollary \ref{cor3.13} it follows that
\begin{equation}\label{4.2}
C \quad = \quad \{ tG_C(u) + (1-t)G_C(-u) : (u,t) \in S \times [0,1] \ \}.
\end{equation}

Now define the continuous function $H_C : S \to E$ by
\begin{equation}\label{4.3}
H_C(u) \quad =_{def} \quad  (G_C(u) \ + \ G_C(-u))/2 \quad = \quad G_C(u) - (r/2)u.
\end{equation}
Hence,
\begin{equation}\label{4.4}
H_C(-u) \quad = \quad H_C(u).    \hspace{3cm}
\end{equation}
That is, $H_C$ is an \emph{even function}.
Furthermore,
\begin{equation}\label{4.5}
G_C(\pm u) \quad = \quad H_C(u) \ \pm \ (r/2)u.   \hspace{3cm}
\end{equation}

Bayen, Lachand-Robert and Oudet call function $H_C$ the \emph{median surface function}. Following them, it will
 be the focus of our parametrization efforts.  Imagine a stick of length
$r$ whose endpoints are the antipodal points. As the stick changes its position in $C$ to point in
direction $u$ the midpoint is at $H_C(u)$.

\begin{theo}\label{theo4.1} If $C$ is a $r$-maximal subset of $E$ then for all $u,v \in S$
\begin{equation}\label{4.6}
(r/2)u \cdot u \ + \ u \cdot H_C(u) \quad \geq \quad (r/2)u \cdot v \ + \ u \cdot H_C(v).
\end{equation}
or, equivalently,
\begin{equation}\label{4.6a}
 u \cdot [H_C(v) - H_C(u)]  \quad \leq \quad (r/4)| u - v |^2.
\end{equation}
\end{theo}

{\bfseries Proof:}  This equation is equivalent to
\begin{equation} \label{4.7}
u \cdot G_C(u) \quad \geq \quad u \cdot G_C(v).
\end{equation}
Because $G_C$ is a surjective function this just says that $\omega_u$ achieves its maximum on
$Bdry(C)$ at $x_1 = G_C(u)$.  This follows from the definition (\ref{3.3}) of $G_C$.

For the second equation we observe that $u  \cdot u  - u \cdot v = (1/2)|u - v|^2$ since $u,v \in S$.

$\Box$ \vspace{.5cm}

\begin{df}\label{def4.2} Let $H : S \to E$ be a continuous, even function.
We say that $H$ satisfies the \emph{$r$-Median Inequality} when for all $u,v \in S$
\begin{equation}\label{4.8}
(r/2) u \cdot u \ + \ u \cdot H(u) \quad \geq \quad (r/2) u \cdot v \ + \ u \cdot H(v).
\end{equation}
We say that $H$ satisfies the \emph{Strict $r$-Median Inequality} when the inequality is strict
whenever $u \not= v$ in $S$.
\end{df}
\vspace{.5cm}

\begin{prop} \label{prop4.2a}  (a) The constant function $H = 0$ satisfies the Strict
$r$-Median Inequality.

(b) Let $H_1, H_2 : S \to E$ be  continuous, even functions with $H_1$ satisfying the $r$-Median Inequality.
If $H_2$ satisfies the $r$-Median Inequality (or the Strict $r$-Median Inequality) then for $0 \leq \lambda < 1$
$\lambda H_1 + (1 - \lambda)H_2$ is a continuous, even function satisfying the $r$-Median Inequality (resp. the
Strict $r$-Median Inequality).
\end{prop}

{\bfseries Proof:} If $u \not= v$ in $S$ then $|u - v|^2 > 0$.  Hence, the zero function satisfies the
Strict $r$-Median Inequality. The second part is obvious.

$\Box$ \vspace{.5cm}

\begin{theo}\label{theo4.3}  Assume that $H : S \to E$ is a continuous, even function.
Define $C \subset E$ by
\begin{equation}\label{4.9}
\begin{split}
C \quad =_{def} \quad \{  H(u) + t(r/2)u : (u,t) \in S \times [-1,1] \ \} \\
= \quad \{   H(u) + t(r/2)u : (u,t) \in S \times [0,1] \ \}. \hspace{1cm}
\end{split}
\end{equation}
If $H$ satisfies the $r$-Median Inequality then $C$ is a $r$-maximal subset of $E$ with $H_C = H$.

If $H$ satisfies the Strict $r$-Median Inequality then, in addition,  $G_C : S \to Bdry(C)$ is a homeomorphism and
each point $x \in Bdry(C)$ has a unique antipodal point in $Bdry(C)$.  In fact, $A_C(x) =
G_C(-(G_C)^{-1}(x))$. Furthermore, at $x = G_C(u)$ there is a unique hyperplane of support for $C$ which contains
$x$ and this hyperplane is perpendicular to $u$.
\end{theo}

{\bfseries Proof:}  The two descriptions of $C$ agree because $H$ is even.

Define $G(u) = H(u) + (r/2)u$. Thus, the $r$-Median Inequality says  $u \cdot G(u) \geq u \cdot G(v)$.

We begin by assuming that $H$ satisfies the Strict $r$-Median Inequality.

 If $v \not= u$ then $u \cdot G(u) > u \cdot G(v)$.  Hence, $G(u) \not= G(v)$.
That is, $G$ is a homeomorphism of the sphere $S$ onto its image $G(S) \subset E$. By the $n-1$ dimensional
Jordan Curve Theorem, $G(S)$ is the boundary of the unique bounded component $U$ of the open set $E \setminus G(S)$.
Let $D = \overline{U} = U \cup G(S)$. On $G(S)$ the functional $\omega_u$ achieves its maximum at $G(u)$ by the
$r$-Median Inequality. The open half-space $\{ y \in E : u \cdot y > u \cdot J(u)\}$ is connected and disjoint
from $G(S)$. Hence, it is disjoint from $D$. That is, $u \cdot G(u)$ is the maximum value of $\omega_u$ on $D$.
It now follows from Proposition \ref{prop2.12} that $D$ is convex. Furthermore, $-u \cdot G(-u)$ is the
maximum value of $\omega_{-u} = - \omega_u$ on $D$. Since,
$G(u) - G(-u) = ru$ it follows that $diam_{\omega_u}(D) = r$.  As $u \in S$ was arbitrary it follows that
$D$ is a closed, convex set of constant diameter $r$ and so $D$ is $r$-maximal.
By definition of $G_D$ we have $G_D(u) = G(u)$.  Hence, $H_D(u) = (G_D(u) + G_D(-u))/2 = (G(u) + G(-u))/2 = H(u)$.
In addition, $C = D$ by equation (\ref{4.2}) because $D$ is $r$-maximal with boundary $G(S)$.  For each $x$ in the
boundary there is a unique $u \in S$ such that $G_C(u) = x$ and so its unique antipode is $G_C(-u)$. Furthermore,
$u = v$ is the unique point $v \in S$ such that $\omega_v$ takes its maximum on $C$ at $x$. Hence, there is a
unique hyperplane of support at $x$ and it is perpendicular to $u$.
This completes
the proof in the strict case.

Now assume that $H$ satisfies $r$-Median Inequality and that $0 \leq \lambda < 1$.

By Proposition \ref{prop4.2a}  $\lambda H = \lambda H + (1 - \lambda) 0$ satisfies the Strict $r$-Median Inequality.
Let $C_{\lambda}$ denote the $2$-round subset
defined by (\ref{4.9}) with $H$ replaced by $\lambda H$.

Now we use the Hausdorff metric on compact subsets of $E$
induced by the metric on $E$. As $\lambda$ approaches $1$ the compact sets $C_{\lambda}$ approach $C$. It follows from
Theorem \ref{theo2.13} that $C$ has constant diameter $r$.  The proof that $G = G_C$ follows by letting $\lambda$
approach $1$ and then the proof that $H_C = H$ follows as before.

$\Box$ \vspace{.5cm}

In the process we have proved the following.

\begin{theo}\label{theo4.4}  Assume that $H : S \to E$ is a continuous, even function which satisfies the
$r$-Median Inequality. For $\lambda \in [0,1]$ define $C_{\lambda} \subset E$ by
\begin{equation}\label{4.11a}
C_{\lambda} \quad =_{def} \quad \{  \lambda H(u) + t(r/2)u : (u,t) \in S \times [-1,1] \ \}.
\end{equation}
For $\lambda \in [0,1)$ the function $\lambda H$ satisfies the Strict $r$-Median Inequality.  As
$\lambda $ varies the $r$-maximal sets $C_{\lambda}$ vary continuously with respect to the Hausdorff metric on
the space of compacta in $E$. With $\lambda = 0$, $C_{\lambda}$ is the $r/2$ ball centered at the origin in $E$. \end{theo}

\vspace{.5cm}
{\bfseries Remark:} Notice that the set of functions $H : S \to E$ which are even and satisfy the $r$-Median Inequality
is closed under pointwise limits.

\vspace{.5cm}

In order to obtain explicit examples where the $r$-Median Inequality holds we introduce smoothness conditions on $H$.
With $G(u) = (r/2)u + H(u)$ defined on $S$ the $r$-Median Inequality says that the function $u \cdot G(v)$ has a local
maximum at $v = u$. Now assume that $H = (h_1,...,h_n)$ is at least $C^1$. Precomposing $H$ and $G$ by the retraction
$r : {\mathbb R}^n \setminus 0 \to S$ given by $r(x) = x/|x|$ we extend $H$ and $G$ and each of their components to
become $C^1$ homogeneous functions of degree zero defined on ${\mathbb R}^n \setminus 0$.
We will denote $\omega_u \circ G$ by $N_u$ so that
\begin{equation}\label{4.11}
N_u(x) \ = \  u \cdot[ (r/2)x /|x|  + H(x) ]  \ = \ (r/2)|x|^{-1} (u \cdot x)  + \sum_i u_i h_i(x).
\end{equation}
 \begin{equation} \label{4.12}
 \frac{\partial N_u}{\partial x_j} \quad = \quad (r/2)[-|x|^{-3} (u \cdot x)x_j  + |x|^{-1}u_j] +
 \sum_i u_i \frac{\partial h_i}{\partial x_j}.
\end{equation}
 The $r$-Median Inequality implies that $N_u = \omega_u \circ G$ has a local maximum at $x = u$.
 Because $x = u$ is a critical point  these partials must vanish.
This says that the image of the tangent map
of $H$ at $u$ is contained in the subspace of vectors perpendicular to $u$, which we will denote $u^{\perp}$.
That is, the image lies in the tangent hyperplane of the sphere $S$ at $u$.
Since this is true for every $u \in S$ and $H$ is homogeneous of degree zero we obtain on ${\mathbb R}^n \setminus 0$
\begin{equation}\label{4.13}
\sum_{i=1}^n \  x_i \cdot \frac{\partial h_i}{\partial x_j} \quad = \quad 0 \qquad \mbox{for} \ j = 1,...,n.
\end{equation}

Now define the real-valued function which is $C^1$ and homogeneous of degree 1:
\begin{equation}\label{4.14}
g(x) \quad =_{def} \quad \sum_{i=1}^n \ x_i h_i(x). \hspace{2cm}
\end{equation}
The above equation is equivalent to:
\begin{equation}\label{4.15}
\frac{\partial g}{\partial x_j} \quad = \quad h_j \qquad \mbox{for} \ j = 1,...,n.
\end{equation}
That is, $H$ is the gradient of $g$. Since the gradient of $g$ is $C^1$ it follows that $g$ is $C^2$.

From (\ref{4.15}) we see that
\begin{equation}\label{4.16}
g(x) \quad =_{def} \quad \sum_{i=1}^n \ x_i \frac{\partial g}{\partial x_i}. \hspace{2cm}
\end{equation}

We recall some well-known homogeneity results.

\begin{lem}\label{lem4.5} Let $g : {\mathbb R}^n \setminus 0 \to {\mathbb R}$ be differentiable.
$g$ is homogeneous of degree $k$ iff
\begin{equation}\label{4.17a}
kg(x) \quad =_{def} \quad \sum_{i=1}^n \ x_i \frac{\partial g}{\partial x_i}. \hspace{2cm}
\end{equation}
In that case, if $g$ is odd (i.e. $g(-x) = -g(x)$) then the gradient is even. Conversely, if the gradient is
even and $k \not= 0$ then $g$ is odd.
\end{lem}

{\bfseries Proof:}  Let $b(t) =_{def} g(tx)$ for $t > 0$.  Then
$b'(t) =  \sum_{i=1}^n \ x_i \frac{\partial g}{\partial x_i}(tx)$.

If $g$ is homogeneous of degree $k$ then $b(t) = t^kb(x)$ and so $b'(t) = kt^{k-1}g(x)$.
Set $t = 1$.

Conversely, if the equation holds then $tb'(t) = kb(t)$ and so $b(t) = t^k \cdot C$ for some constant $C$ and setting
$t = 1$ we see that the constant $C = g(x)$.

If $g(-x) = -g(x)$ then taking partials we get that each $\frac{\partial g}{\partial x_i}$ is even.
The converse follows from the
homogeneity equation (\ref{4.17a}) when $k \not= 0$.

$\Box$ \vspace{.5cm}

Thus, a $C^1$  function $H$ which is homogeneous of degree zero and whose tangent map at $x$
has image in $x^{\perp}$ for
every $x$ is exactly the gradient of a $C^2$ real-valued function $g$ which is odd and homogeneous of degree one,
i.e. which satisfies $g(tx) = tg(x)$ for all $t \not= 0$.  Extending continuously, but not smoothly, to the origin
we have $g(tx) = tg(x)$ for all real $t$.

Thus, equation (\ref{4.12}) becomes
 \begin{equation} \label{4.17}
 \frac{\partial N_u}{\partial x_j} \quad = \quad (r/2)[-|x|^{-3} (u \cdot x)x_j  + |x|^{-1}u_j] +
 \sum_i u_i \frac{\partial^2 g}{\partial x_i \partial x_j}.
\end{equation}
At $x = u$ the equation $ \sum_i u_i \frac{\partial^2 g}{\partial x_i \partial x_j} = 0$ holds  because
$\frac{\partial g}{ \partial x_j}$ is homogeneous of degree zero.

\begin{df}\label{def4.6} Let $g : {\mathbb R}^n \setminus 0 \to {\mathbb R}$ be an odd, $C^2$
function which is homogeneous of degree $k = 1$.  We say that $g$ satisfies the \emph{Linear $r$-Median Condition}
if at every point $u \in S$ the eigenvalues of the  Hessian matrix
$(  \frac{\partial^2 g}{ \partial x_k \partial x_j} )$ (hereafter denoted $\partial^2g$)
are contained in the closed interval $[-(r/2),(r/2)]$.
We say that  $g$ satisfies the \emph{Strict $r$-Median Condition}
if at every point $u \in S$ the eigenvalues of $\partial^2g$ are contained in the open interval $(-(r/2),(r/2))$.
\end{df} \vspace{.5cm}

{\bfseries Remark:}(a) The Hessian matrix function $\partial^2g$ on ${\mathbb R}^n \setminus 0$ is homogeneous of
degree $-1$. This is why the eigenvalue conditions are restricted to points $u$ on the sphere.  The eigenvalues
will blow up as $x$ approaches the origin.

(b) The function $\partial^2g$ is odd and so if $ \mu$ is an eigenvalue at $u$ then $- \mu $ is an eigenvalue at $-u$.
\vspace{.5cm}

\begin{prop}\label{prop4.6a} Let $g :S \to {\mathbb R}$ be any odd, $C^2$ function on the unit sphere. Extend
$g$ to an odd, $C^2$ homogeneous function of degree $1$ on ${\mathbb R}^n \setminus 0$ by $g(x) = |x| g(x/|x|)$.
There is a maximum positive $\lambda^*$ such that $\lambda^* g$ satisfies the  Linear $r$-Median Condition and
then $\lambda g$ satisfies the Strict Linear $r$-Median Condition for all $\lambda \in [0,\lambda^*)$.\end{prop}

{\bfseries Proof:}  Multiplying $g$ by $\lambda$  multiplies the eigenvalues of the Hessian at $u$ by $\lambda$
for every $u \in S$. As the eigenvalues vary continuously the result is obvious from compactness of $S$.

$\Box$ \vspace{.5cm}

\begin{lem} \label{lem4.6b} For $g : {\mathbb R}^n \setminus 0 \to {\mathbb R}$  an odd, $C^2$
function which is homogeneous of degree $1$ let $H : {\mathbb R}^n \setminus 0 \to {\mathbb R}^n$ be the
gradient of $g$.  If $H$ satisfies the $r$-Median Inequality, then $g$ satisfies the Linear $r$-Median Condition.
\end{lem}

{\bfseries Proof:}  For any point $u \in S$ the vector $u$ itself is an eigenvector of
$\partial^2g$ at $u$ with eigenvalue zero.
The eigenvectors with nonzero eigenvalues lie in $u^{\perp}$.  Let $v \in u^{\perp}$ and
define the path $p : [0,\pi] \to S$ by
\begin{equation}\label{4.17aa}
\begin{split}
p(t) \quad = \quad  \cos(t) u  \ + \ \sin(t) v, \quad \mbox{so that} \\
p\ '(t) \quad = \quad - \sin(t) u \ + \ \cos(t) v \hspace{1cm}
\end{split}
\end{equation}
At $p(t)$, the vector $p(t)$ is an eigenvector of $\partial^2g$ with eigenvalue zero and so at $p(t)$
\begin{equation}\label{4.18a}
 \sum_i u_i \frac{\partial^2 g}{\partial x_i \partial x_j} \quad =
 \quad - \tan(t)  \sum_i v_i \frac{\partial^2 g}{\partial x_i \partial x_j}.
\end{equation}
Differentiating $N_u(p(t))$ we get, after applying (\ref{4.18a}) on each side of $\partial^2 g$
\begin{equation}\label{4.19a}
\begin{split}
\frac{dN_u \circ p}{dt} \ = \ (r/2) u \cdot p\ '(t) \
+ \  \sum_{i,j} u_i \frac{\partial^2 g}{\partial x_i \partial x_j}p\ '(t)_j \ = \\
- (r/2) \sin(t) \ - \ \tan(t)  \sum_{i,j} v_i \frac{\partial^2 g}
{\partial x_i \partial x_j}[\sin(t) \tan(t) + \cos(t)]v_j \\
= - \sin(t)[ (r/2) + \sec(t)(\sin(t) \tan(t) + \cos(t))
\sum_{i,j} v_i \frac{\partial^2 g}{\partial x_i \partial x_j}v_j].
\end{split}
\end{equation}

By the $r$-Median Inequality $N_u$ has a local maximum at $u$. So by the Mean Value Theorem there
is a positive sequence $\{ t_k \}$ tending
to zero such that $d(N_u \circ p)/dt \leq 0$ at each point $p(t_k)$.   Applying (\ref{4.19a}),
dividing by the positive number
$\sin(t_k)$ and then letting $k$ tend to infinity we obtain at $u = p(0)$:
\begin{equation}\label{4.20a}
(r/2) + \sum_{i,j} v_i \frac{\partial^2 g}{\partial x_i \partial x_j}v_j \quad \geq \quad 0
\end{equation}
for all $v \in u^{\perp}$.  In particular,  every nonzero eigenvalue $\mu$ of $\partial^2g$
at $u$ satisfies $\mu \geq - (r/2)$.
Since $- \mu$ is an eigenvalue at $-u$ we have $(r/2) \geq \mu \geq -(r/2)$ as required.

$\Box$ \vspace{.5cm}

\begin{theo}\label{theo4.7a}  For $g : {\mathbb R}^n \setminus 0 \to {\mathbb R}$  an odd, $C^2$
function which is homogeneous of degree $1$ let $H : {\mathbb R}^n \setminus 0 \to {\mathbb R}^n$ be the even, $C^1$
function which is the gradient of $g$ so that $H$ is homogeneous of degree $0$.

 The function $g$ satisfies the Strict Linear $r$-Median Condition iff $-r/2$ is not an eigenvalue of the Hessian matrix
 $(  \frac{\partial^2 g}{ \partial x_k \partial x_j} )$ at any point $u$ of $S$.

 If $g$ satisfies the Strict Linear $r$-Median Condition then $H$ satisfies the Strict $r$-Median Inequality.

 Furthermore, $g$ satisfies the Linear $r$-Median Condition iff $H$ satisfies the $r$-Median Inequality.\end{theo}

 {\bfseries Proof:} It is obvious that if $g$ satisfies the Strict Linear $r$-Median Condition
 then $-r/2$ is not an eigenvalue of $\partial^2g$ at any point of $S$.

  Now assume that $-r/2$ is not an eigenvalue of $\partial^2g$ at any point
 of $S$.

Let $\bar G = (r/2)I + H$ which agrees with $G$ on $S$. For $\bar G$ we are not projecting to
 the sphere for the $r/2$ times the identity term and so $\bar G $ is not homogeneous. The assumption that
 $-r/2$ is never an eigenvalue says exactly that the tangent map of $\bar G$ is injective at each point $v$ of the
 sphere. Meanwhile, on the sphere the tangent map of $G$ at $v$ maps $v^{\perp}$ into itself.  Since this is the
 restriction of the tangent map of $\bar G$ as well, it follows that at each point $v \in S$ the tangent map
 of $G$ is an isomorphism of $v^{\perp}$.

 Now consider  $N_u = \omega_u \circ G$ on $S$.  At any point $v \not= \pm u$, $\omega_u$ maps $v^{\perp}$ onto
 ${\mathbb R}$.  Hence, $\pm u$ are the only critical points of $\N_u$ on $S$.  On the other hand, the
 maximum and minimum values of $N_u$ on $S$ are certainly critical values. Since
 $\omega_u(G(u)) - \omega_u(G(-u)) = r > 0$ it follows that $N_u(u)$ and $N_u(-u)$ are the maximum
 and minimum values of $N_u$ on $S$. Furthermore, no other points are critical and so $N_u$
 achieves its maximum only at $u$.  This proves that $H$ satisfies the Strict $r$-Median Inequality.

 To complete the proof of the first part, we must show that $g$ satisfies the Strict Linear $r$-Median Condition.

Let $\mu$ is an eigenvalue of $\partial^2g$ at some point.
 In that case, Lemma \ref{lem4.6b} implies that $\mu$ lies in the interval $[-(r/2),(r/2)]$.  By assumption
 $\mu $ is not $-(r/2)$. As $- \mu$ is  an eigenvalue at the antipodal point, $ - \mu$ is not $-(r/2)$ either.
 Hence, $\mu$ is in the open interval and the Strict Linear $r$-Median Condition holds.

Now if $g$ satisfies the Linear $r$-Median Condition then for $\lambda < 1$, $\lambda g$ satisfies the Strict
Linear $r$-Median Condition and so $\lambda H$ satisfies the Strict $r$-Median Inequality.  Letting $\lambda$ approach
$1$ we see that $H$ satisfies the $r$-Median Inequality.

Conversely, if $H$ satisfies the $r$-Median Inequality then $g$ satisfies the Linear $r$-Median
Condition by Lemma \ref{lem4.6b} again.

$\Box$ \vspace{.5cm}

\begin{theo}\label{theo4.7b} Let $g : {\mathbb R}^n \setminus 0 \to {\mathbb R}$  an odd, $C^2$
function which is homogeneous of degree $1$.  There exists a nonnegative real number $r^*$ such that
\begin{equation}\label{4.21a}
\begin{split}
\{ \mu : \mu \mbox{\ is an eigenvalue of the matrix \ } \partial^2g \\
\mbox{ \ at some \ } u \in S \ \} \quad = \quad
[-(r^*/2),(r^*/2)].
\end{split}\end{equation}

Let $H$ be the gradient of $g$ and let $r$ be a positive real number.  If $r > r^*$ then $H$ satisfies the Strict
$r$-Median Inequality.  If $r = r^*$ then $H$
satisfies the $r$-Median Inequality.  If $r < r^*$  then $H$ does not satisfy the
$r$-Median Inequality.
\end{theo}

{\bfseries Proof:} Since the eigenvalues of $\partial^2g$ vary continuously on the compact set $S$ it follows
that the set of eigenvalues is compact as $u$ varies over $S$. At every point $u \in S$, $u$ is an eigenvector
of $\partial^2g$ with eigenvalue $0$. Furthermore, if $\mu$ is an eigenvalue at $u$ then $- \mu$ is an
eigenvalue at $-u$.  Hence, there is a nonnegative $r^*$ such that the set of eigenvalues is contained in
the closed interval $[-(r^*/2),(r^*/2)]$ and the endpoints, $\pm r^*/2$ and also $0$ are eigenvalues. Now
suppose $r > 0$ with $- r/2$ not an eigenvalue.   Theorem \ref{theo4.7a} implies that
the set of eigenvalues is contained in $(-(r/2),(r/2))$ and so $r > r^*$. Moreover, if $r/2$ is not an eigenvalue
then by symmetry $- r/2$ is not an eigenvalue.  Thus, the set of eigenvalues is the entire closed interval
$[-(r^*/2),(r^*/2)]$.

Clearly, $g$ satisfies the Linear $r$-Median Condition iff $r \geq r^*$ and the Strict Linear $r$-Median Condition
iff $r > r^*$.  Thus, the remaining results follow from Theorem \ref{theo4.7a}.

$\Box$ \vspace{.5cm}

Thus, any $C^2$ odd, homogeneous $g$ satisfies the Linear $r$-Median Condition for $r$ sufficiently large.
We can use $r$
as a parameter for the solid $C$ described by (\ref{4.9}), with $H$ the gradient of $g$.
 At $r= 0$, $C$ is
just the image  $H(S)$ and for small $r \geq 0 $  the  Linear $r$-Median Condition will not hold.  The set $C$ will not
be convex. For $r \geq r^*$ sufficiently large, $C$ becomes convex and the Linear $r$-Median Condition holds.  So for all
$r \geq r^*$ the set $C$ is a solid of constant diameter $r$.

Notice that the trivial case of an odd function $g$ of degree $1$ is the linear function $ a \cdot x$ for
a constant vector $a$. Adding this linear function to $g$ translates the gradient $H$ by  $a$ and leaves the
Hessian matrices unaffected.  In particular the Linear $r$-Median Conditions are unaffected. The resulting $r$-maximal
sets are translated in ${\mathbb R}^n$.

If $O$ is an orthogonal $n \times n$ matrix, let $u \mapsto O(u)$ denote the associated linear map.
If $g_O = g \circ O^{-1}$ then the gradient $H_O$ satisfies $H_O(u) = O(H(O^{-1}(u)))$.  Hence,
\begin{equation}\label{4.21b}
\begin{split}
G_O(u) \ = \ (r/2)u + H_O(u) \ = \hspace{1cm}\\ O((r/2)O^{-1}(u) + H(O^{-1}(u)) \ = \ O(G(O^{-1}(u))).
\end{split}\end{equation}
\vspace{.5cm}
Thus, the image of $G_O$ is the image of $G$ rotated by $O$.
\vspace{.5cm}

In the case of the Euclidean Plane, i.e. $ n = 2$, there is a somewhat different parametric approach.

We use polar coordinates $\rho, \t$ for ${\mathbb R}^2 \setminus 0$. Thus, the angular coordinate $\t$
parametrizes the unit circle $S$. A function $H = (h_1,h_2) : {\mathbb R}^2 \setminus 0 \to {\mathbb R}^2$
is homogeneous of degree $0$ when it is
a function of $\t$ alone.  A function $g : {\mathbb R}^2 \setminus 0 \to {\mathbb R}$
 is homogeneous of degree $1$ when there exists a function
$a : {\mathbb R}^2 \setminus 0 \to {\mathbb R}$ homogeneous of degree $0$ with $g(\rho,\t) = \rho a(\t)$.
$H$ is an even function when $a(\t + \pi) = a(\t)$ for all $\t$,i.e. $a$ has period $\pi$ and $H$ is an odd function when
$a(\t + \pi) = - a(\t)$, i.e. $a$ has \emph{anti-period} equal to $\pi$.
For example, $\cos(k \t)$ and $\sin(k \t)$ are both $\pi$ periodic when $k$ is even and $\pi$ anti-periodic
when $k$ is  odd.

With $x = \rho \cos(\t), y = \rho \sin(\t)$ the coordinate change matrices are given by
\begin{equation} \label{4.22}
 \begin{pmatrix} \partial x/\partial \rho & \partial x/\partial \t\\
\partial y/\partial \rho & \partial x/\partial \t \end{pmatrix}  \quad = \quad
\begin{pmatrix} \cos(\t) & -\rho \sin(\t) \\ \sin(\t) & \rho \cos(\t) \end{pmatrix} ,
\end{equation}
and its inverse
\begin{equation} \label{4.23}
 \begin{pmatrix} \partial \rho/\partial x & \partial \rho/\partial y\\
\partial \t/\partial x & \partial \t/\partial y \end{pmatrix}  \quad = \quad
\begin{pmatrix} \cos(\t) & \sin(\t) \\ - \rho^{-1}\sin(\t) & \rho^{-1} \cos(\t) \end{pmatrix} ,
\end{equation}

Now consider $g(\rho,\t) = \rho a(\t)$ with $a$ a $C^2$ function of the angular coordinate.
The gradient $H = (h_1,h_2)$ is given by:
\begin{equation}\label{4.24}
\begin{split}
h_1(\t) = \frac{\partial g}{\partial \t} \frac{\partial \t}{\partial x} +
\frac{\partial g}{\partial \rho}\frac{\partial \rho}{\partial x} = - a'(\t) \sin(\t) + a(\t) \cos(\t), \\  \\
h_2(\t) = \frac{\partial g}{\partial \t} \frac{\partial \t}{\partial y} +
\frac{\partial g}{\partial \rho}\frac{\partial r}{\partial y} =  \ a'(\t) \cos(\t) + a(\t) \sin(\t).
\end{split}
\end{equation}

Now define $\b$ so that
\begin{equation}\label{4.25}
(r/2) \b(\t) \quad =_{def} \quad a(\t) \ + \ a''(\t). \hspace{2cm}
\end{equation}

Regarding the restriction of $H$ to $S$ as a path in ${\mathbb R}^2$ we see that its velocity vector is
given by
\begin{equation}\label{4.26}
H'(\t) \quad = \quad (h_1'(\t), h_2'(\t)) \quad = \quad (r/2) \b(\t) (- \sin(\t), \cos(\t)).
\end{equation}

Finally, the Hessian matrix for $g$ is given by
\begin{equation}\label{4.27}
\begin{split}
 \begin{pmatrix} \partial^2 g/\partial x^2 & \partial^2 g/\partial x\partial y\\
\partial^2 g/\partial y \partial x & \partial^2 g/\partial y^2 \end{pmatrix}  \quad = \quad
 \begin{pmatrix} \partial h_1/\partial x & \partial h_1/\partial y\\
\partial h_2/\partial x & \partial h_2/\partial y \end{pmatrix}  \\  \\
= \quad  \rho^{-1}(r/2)\b(\t)  \begin{pmatrix} \sin^2(\t) &  - \sin(\t) \cos(\t) \\
-\sin(\t) \cos(\t)  &  \cos^2(\t)  \end{pmatrix} \hspace{2cm}
\end{split}
\end{equation}

The two eigenvalues of the latter matrix are $0$ and $\rho^{-1}(r/2) \b(\t)$.  From all this we obtain

\begin{prop}\label{prop4.8}  Assume that $a : S \to {\mathbb R}$ is a $\pi$ anti-periodic, $C^{k+2}$ function. The
function $ \b : S \to {\mathbb R}$ defined by $\b(\t) = (2/r)[a(\t) + a''(\t)]$
is a $\pi$ anti-periodic, $C^k$ function
such that
\begin{equation}\label{4.28}
\int_0^{\pi} \ \b(\t)(- \sin(\t), \cos(\t)) d \t \quad = \quad (0,0).
\end{equation}
Furthermore, the odd function, homogeneous of degree $1$, defined by $g(\rho,\t) = \rho a(\t)$ satisfies the
Linear $r$-Median Condition (or the Strict Linear $r$-Median Condition) iff $|\b(\t)| \leq 1$ for all $\t \in {\mathbb R}$
(resp. $|\b(\t)| < 1$ for all $\t \in {\mathbb R}$ ).

Conversely, given a $\pi$ anti-periodic $C^k$ function  $\b(\t)$ which satisfies (\ref{4.28}) there exists
a $\pi$ anti-periodic, odd $C^{k+1}$ function $a(\t)$ such that $(r/2) \b = a + a''$.
\end{prop}

{\bfseries Proof:} Equation (\ref{4.28}) follows from (\ref{4.26}) because integrating $H'$ from $0$ to $\pi$
 yields $H( \pi) - H(0)$ which equals   $ 0$ because $H$ is even. It is clear that if $a$ is $\pi$ anti-periodic then
 $\b = (2/r)[a + a'']$ is $\pi$ anti-periodic. The Linear $r$-Median Conditions
 follow from (\ref{4.27}) because with $r = 1$ the only
 nonzero eigenvalue at $\t$ is $(r/2) \b(\t)$.

 Conversely, given $\b$ we can use variation of parameters to solve the differential equation:
 \begin{equation}\label{4.29}
 a(\t) \quad =_{def} \quad \int_0^{\t} \  (r/2) \b(s)[\cos(s)\sin(\t) - \sin(s) \cos(\t)] \ ds.
\end{equation}
By (\ref{4.28}) $a(0) = 0 = a(\pi)$.  Because $\sin(\t)$ and $\cos(\t)$ are $\pi$ anti-periodic it follows that
 \begin{equation}\label{4.30}
 a(\t + \pi) \quad = \quad - \int_0^{\t + \pi} \ (r/2) \b(s)[\cos(s)\sin(\t) - \sin(s) \cos(\t)] \ ds.
\end{equation}
On the other hand, $\b(s)(-\sin(s),\cos(s))$ is $\pi$ periodic. Hence, in (\ref{4.28})
the integral from $\t$ to $\t + \pi$ is zero for any $\t$.  It follows that $a(\t + \pi) = - a(\t)$.
Thus, $a$ is $\pi$ anti-periodic.

The general solution of the differential equation is obtained by adding $C_1 \cos(\t) + C_2 \sin(\t)$ for
arbitrary constants $C_1, C_2$.  Adding such a term corresponds to adding a linear function to $g(r,\t)$
and has the effect of translating the gradient $H$ by a constant.

$\Box$ \vspace{.5cm}

Above we used $g$ or, equivalently, $a(\t)$ to construct our parametrization.  Now we will use $\b(\t)$
instead. As we are dealing with curves in the plane we will switch to the usual parametric representation for
such curves.  For example, $t \mapsto \mathbf{U}(t) =_{def} (\cos (t),\sin (t))$ parametrizes the unit circle
with unit tangent vector $\mathbf{T}(t) =_{def} (- \sin (t),\cos (t))$.
We also return now to the general case of $r$-rotundity.

We begin with a $C^2$ function $\b : {\mathbb R} \to {\mathbb R}$ such that
\begin{itemize}
\item[(i)]$|\b(t)| \quad \leq \quad 1.$
\item[(ii)] $\b(t + \pi) \quad = \quad - \b(t).$
\item[(iii)] $\int_0^{\pi} \ \b(u) ( - \sin(u),\cos(u)) \ du \quad = \quad (0,0).$
\end{itemize}

The associated boundary curve (the one dimensional version of  \\ $G(u)$ above)
has position and velocity vectors given by
\begin{equation}\label{4.31}
\begin{split}
\mathbf{R}(t) \quad =_{def} \quad  \int_0^t \ (r/2)(1 + \b(u))(- \sin(u), \cos(u)) \ du \\
\mathbf{R}'(t) \quad = \quad (r/2) (1 + \b(t))(- \sin(t), \cos(t)). \hspace{1cm}
\end{split}
\end{equation}
Thus, the unit tangent vector is $\mathbf{T}(t)$ and the speed is given by
\begin{equation}\label{4.32}
\frac{ds}{dt} \quad = \quad (r/2)(1 + \b(t)).
\end{equation}
Notice that $\mathbf{T}'(t) = - \mathbf{U}(t)$ and so $- \mathbf{U}(t)$ is the unit normal vector.
Hence,
\begin{equation}\label{4.33}
\frac{d\mathbf{T}}{ds} \quad = \quad - [(r/2)(1 + \b(t))]^{-1} \mathbf{U}(t)
\end{equation}

By Proposition \ref{prop4.8}, Theorem \ref{theo4.7a}, and Theorem \ref{theo4.4} the closed curve
bounds the $r$-maximal set in the plane:
\begin{equation}\label{4.34}
C_{\b} \quad =_{def} \quad \{ u \mathbf{R}(t) + (1-u) \mathbf{R}(t + \pi) : (t,u) \in [0, 2 \pi] \times [0,1] \ \}
\end{equation}

Notice that for the arclength of the closed curve we  get
\begin{equation}\label{4.34a}
\int_0^{2\pi} \ (r/2)(1 + \b(t)) \ dt \quad = \quad \pi r.
\end{equation}
because the function $\b$ is $\pi$ anti-periodic.  This is the smooth special case
of a general theorem due to Barbier, for its
proof, see Lyusternik (1966) Section 12.

\begin{theo}\label{theo4.9} Let $\b : [0, \pi) \to [-1,1]$ be a measurable function such that
$\int_0^{\pi} \ \b(u) ( - \sin(u),\cos(u)) \ du \ =  \ (0,0)$. Extend $\b$ to $[0,2\pi)$ by $\b(t + \pi) =
- \b(t)$ for $t \in [0,\pi)$. Then extend to ${\mathbb R}$ to make the function $2 \pi$ periodic.
The Fourier series has nonzero coefficients only for
$\sin (kt)$ and $\cos (kt)$ with $k$ odd and greater than $1$.

The curve $\mathbf{R}(t) = \int_0^t \ (r/2)(1 + \b(u))(- \sin(u), \cos(u)) \ du $ is continuous and
bounds an $r$-maximal subset $C_{\b}$.
\end{theo}

{\bfseries Proof:} One can construct a sequence of smooth functions $\b_n$ which satisfy conditions
(i), (ii) and (iii) and which
converge to $\b$ in $L^1([0,2\pi])$. The characteristic functions of $[0,t]$, $\sin (t)$ and $\cos (t)$ all
have norm at most 1 in $L^{\infty}$. Hence, the curves for $\b_n$ converge uniformly to the curve for
$\b$. The latter is continuous but need not be differentiable. Each $C_{\b_n}$  is an $r$-maximal set
and as $n \to \infty$ these sets converge to  $C_{\b}$ with respect to the Hausdorff metric. By Theorem
\ref{theo2.13} $C_{\b}$ is $r$-maximal.

$\Box$ \vspace{.5cm}

For the regular polygon with $2k + 1$ vertices we use the function $\beta(t)$
\begin{equation}\label{4.35}
\b(t) \quad = \quad (-1)^i \qquad \mbox{for} \quad t \in
[\frac{i\pi}{2k+1},\frac{(i+1)\pi }{2k+1}) \quad \mbox{with} \quad i = 0,...2k,
\end{equation}
and extend over $(\pi,2\pi)$ so that condition (ii) holds. Notice that for $t \in [0,\pi/2]$ we have
\begin{equation}\label{4.36}
\b(\pi - t) \quad = \b(t).
\end{equation}
So it is clear that $\int_0^{\pi} \b(t) \cos(t) dt = 0$.  The sine condition  is not obvious.
It is true since the example works.  However, there is a geometric argument.
\begin{equation}\label{4.37}
\int_0^{\pi} \beta(t) {\mathbf T}(t) \ dt \quad = \quad \sum_{i=0}^{2k}
(-1)^i[{\mathbf U}(\frac{(i+1)\pi}{2k+1})
- {\mathbf U}(\frac{i \pi}{2k+1})].
\end{equation}
For $i$ odd we replace $- [{\mathbf U}(\frac{(i+1)\pi}{2k+1})
- {\mathbf U}(\frac{i \pi}{2k+1})]$ by $[{\mathbf U}(\frac{(i+1)\pi}{2k+1} + \pi)
- {\mathbf U}(\frac{i \pi}{2k+1} +\pi)] = [{\mathbf U}(\frac{(2k+1+i+1)\pi}{2k+1})
- {\mathbf U}(\frac{2k+1+i \pi}{2k+1})]$ and so we have
\begin{equation}\label{4.38}
\int_0^{\pi} \beta(t) {\mathbf T}(t) \ dt \quad = \quad \sum_{i=0}^{2k}
[{\mathbf U}(\frac{(2i+1)\pi}{2k+1})
- {\mathbf U}(\frac{2i \pi}{2k+1})].
\end{equation}
This equals zero because  the sum is invariant under rotation by $\frac{2\pi}{2k+1}$.

With $k = 1$ we have the triangle case.

In general, let $\mu$ be the measure on $[0,\pi/2]$ with density $\sin(t)$ so that $\mu([0,\pi/2]) = 1$.
For $\b(t)$ defined on $[0,\pi/2]$ with $|\b(t)| \leq 1$ and
such that $\int_0^{\pi/2} \ \b(t) \ \mu(dt) = 0$ we can extend
to $[\pi/2,\pi)$ so that (\ref{4.36}) holds. Then extend to $(0,2\pi)$ to get condition (ii). For example,
if $A$ is any measurable subset of $[0.\pi/2]$ with $\mu(A) = \frac{1}{2}$ then we can use $\b(t) =
\chi_A(t) - \chi_{A'}(t)$ with $A' = [0,\pi/2] \setminus A$.  The interesting case here is with $A$ a Cantor set.

  Any curve with nonvanishing curvature can be parametrized by using the unit normal.
  This is the one dimensional version  of the classical Gauss map.

  \begin{theo}\label{theo4.10} Given a plane curve with  radius of curvature bounded by $r$
   at every point then any sufficiently
  short piece can be embedded in the boundary of some $r$-maximal subset of the plane.  To be precise, if the curve
  is parametrized by $t = $the angle $\t$ of the unit normal and the parameter moves through an interval of length
  at most $\pi/3$ then the curve can be embedded within the boundary of an $r$-maximal subset. \end{theo}

  {\bfseries Proof:} By rotating we can assume that the angle $t$ moves from $0$ to some $\theta^* \leq \pi/3$.
   We can reverse the construction to
  define the associated function $\b$ on $[0,\theta^*]$.

  Since the curve is parametrized by  the angle of the unit normal, the velocity vector is a multiple of
  ${\mathbf T}(t)$.  We define $\b(t)$ so that $(r/2)(1 + \b(t))$ is the speed at $t$.  As above, the radius of
  curvature is given by $(r/2)( 1 + \b(t))$ and since this is between $0$ and $r$ we have $|\b(t)| \leq 1$.
  Since $\mu([0,\pi/3]) = \frac{1}{2}$, we can extend  $\b$ to $(\theta^*,\pi/2)$ by a suitable
  constant in $[-1.1]$ so that
  $\int_0^{\pi/2} \ \b(t) \mu(dt) = 0$ and proceed as above.

  $\Box$ \vspace{.5cm}

  {\bfseries Remark:}   The $\pi/3$ restriction is needed.  Notice that you can't use an
  arc of a circle with radius $r$ and angle greater than   $\pi/3$ because then
 the chord is longer than the radius.
\vspace{1cm}

\section*{References}

E. Akin (1993) {\bfseries The general topology of dynamical
systems} Amer. Math. Soc., Providence.
\vspace{.5cm}

G. Averkov and H. Martini (2004) \emph{A characterization of constant width in Minkowski planes} Aequationes Math.
{\bfseries 68}:38-45. \vspace{.5cm}

T. Bayen, T. Lachand-Robert and E. Oudet (2007)
\emph{Analytic parametrization of three dimensional bodies of constant width} Arch.
Rational Mech. Anal. {\bfseries 186:}225-249.
\vspace{.5cm}

T. Bonnesen and W. Fenchel (1934) {\bfseries Theorie der konvexen k\"{o}rper} Springer-Verlag, Berlin.
\vspace{.5cm}

T. Bonnesen and W. Fenchel (1987) {\bfseries Theory of convex bodies} (translation of Bonnesen
and Fenchel (1934)) BCS Associates, Moscow, Idaho.
\vspace{.5cm}

H. G. Eggleston (1958) {\bfseries Convexity} Cambridge University Press, London. \vspace{.5cm}

H. G. Eggleston (1965) \emph{Sets of constant width in finite dimensional Banach spaces}
Israel J. Math. {\bfseries 3}:163-172. \vspace{.5cm}

B. Guilfoyle and W. Klingenberg (2009)
\emph{On $C^2$ smooth surfaces of constant width} Tbil. Math. J. {\bfseries 2:} 1 - 17.
\vspace{.5cm}

P. C. Hammer (1963) \emph{Convex curves of constant Minkowski width} in {\bfseries Convexty}, ed. V. L. Klee,
Proc. of Symp. Pure Math. VII :291-304.\vspace{.5cm}

K. Kuratowski (1968) {\bfseries Topology} Academic Press, New York.
\vspace{.5cm}

T. Lachand-Robert and E. Oudet (2007)
\emph{Bodies of constant width in arbitrary dimension} Math.
Nachr. {\bfseries 290} No. 7: 740-750.
\vspace{.5cm}

S. R. Lay (1982) {\bfseries Convex sets and their applications} John Wiley \& Sons, New York. \vspace{.5cm}

L. Lyusternik (1966) {\bfseries Convex figures and polyhedra} D. C. Heath, Chicago.
\vspace{.5cm}

\end{document}